\documentclass[12pt]{amsart}
\usepackage{amsfonts, amssymb,amsmath, amsthm, amsxtra, latexsym, mathrsfs, amscd}
\usepackage{amsmath, amssymb, amsthm, times, float}
\usepackage{mathtools, amssymb, amsthm, times, float}
\usepackage{graphics}
\usepackage{MnSymbol}

\usepackage{amssymb,hyperref}
\usepackage[latin1]{inputenc} %

\usepackage[thinlines]{easytable}
\usepackage{tikz-cd}

\newtheorem{theo+}{Theorem}[section]
\newtheorem{prop+}[theo+]{Proposition}
\newtheorem{coro+}[theo+]{Corollary}
\newtheorem{lemm+} [theo+]{Lemma}
\newtheorem{deep+}  [theo+]  {Deep Result}
\newtheorem{fact+}  [theo+]  {Fact}
\theoremstyle{definition}
\newtheorem{exam+}  [theo+]  {Example}
\newtheorem{rema+}  [theo+]  {Remark}
\newtheorem{defi+}  [theo+]  {Definition}
\newtheorem{xca+}[theo+]{Exercise}

\numberwithin{equation}{section}

\usepackage{color}

\newcommand\beq{\begin{equation}\label}
\newcommand\eeq{\end{equation}}
\renewcommand\({\Big(}
\renewcommand\){\Big)}

\renewcommand\]{\Big]}

\renewcommand\a[1]{{\acute{#1}}}

\def\draft{\centerline{(Draft {\the \day}/{\the\month} \the \year.)}}
\bigskip

\def\refn#1.#2{\expandafter\def\csname#1\endcsname{[#2]}}
\def\refnr#1.{\csname#1\endcsname}

\def\fa{\mathfrak a}

\def\fe{\mathfrak e}

\def\fg{\mathfrak g}
\def\fk{\mathfrak k}
\def\fh{\mathfrak h}
\def\fl{\mathfrak l}
\def\fm{\mathfrak m}
\def\fn{\mathfrak n}
\def\fp{\mathfrak p}
\def\fq{\mathfrak q}

\def\fs{\mathfrak s}
\def\fsl{\mathfrak sl}

\def\fu{\mathfrak u}
\def\fso{\mathfrak{so}}
\def\fsp{\mathfrak{sp}}

\def\fsu{\mathfrak{su}}
\def\fgl{\mathfrak{gl}}

\def\a{\alpha}

\def\Claminv2{|C(\Lambda)|^{-2}}

\def\de{d\varepsilon}

\def\Aa2D{A^{\a,2}(D)}
\def\bAa2D{\overline{A^{\a,2}(D)}}
\def\Ab2D{A^{\beta,2}(D)}
\def\bAb2D{\overline{A^{\beta,2}(D)}}

\def\Norm#1_#2{\Vert#1\Vert_{#2}}

\def\phipl12{\phi_{p_{l_1}, p_{l_2}}}
\def\phip01{\phi_{p_{0}, p_{0}}}

\def\a{\alpha}

\def\Claminv2{|C(\Lambda)|^{-2}}

\def\sig{\sigma}

\def\ad{\operatorname{ad}}
\def\Ad{\operatorname{Ad}}

\def\ch{\operatorname{ch}}

\def\rank{\operatorname{rank}}
\def\diag{\operatorname{diag}}

\def\Ind{\operatorname{Ind}}
\def\diag{\operatorname{diag}}
\def\genus{\operatorname{genus}}

\def\ch{\operatorname{ch}}
\def\exp{\operatorname{exp}}

\def\mod{\operatorname{mod}}

\def\im{\operatorname{Im}}
\def\re{\operatorname{Re}}

\def\sh{\operatorname{sh}}

\def\sh{\operatorname{sh}}

\def\tr{\operatorname{tr}}

\def\de{d\varepsilon}

\def\Aa2D{A^{\a,2}(D)}
\def\bAa2D{\overline{A^{\a,2}(D)}}
\def\Ab2D{A^{\beta,2}(D)}
\def\bAb2D{\overline{A^{\beta,2}(D)}}

\def\phipl12{\phi_{p_{l_1}, p_{l_2}}}
\def\phip01{\phi_{p_{0}, p_{0}}}

\def\alg/{algebra}

\def\Alg/{Algebra}

\def\alt/{alternative} 
\def\anal/{analytic}
\def\analfunc/{\anal/\ \func/}
\def\Ans/{\it Answer. \normal}
\def\ass/{associative}
\def\nass/{non-\ass/}
\def\autom/{automorphism}
\def\homom/{homomorphism}
\def\isom/{isomorphism}
\def\bdd/{bounded}
\def\Bdd/{Bounded}
\def\bddsymdom/{bounded \sym/ \dom/}
\def\Cartdom/{Cartan \dom/}
\def\bdry/{boundary}
\def\bsd/{\bdd/ \symdom/}
\def\bv/{boundary value}
\def\cf/{{\it cf}\.}
\def\Cf/{{\it Cf}\.}
\def\charr/{character}
\def\coeff/{coefficient}
\def\comm/{commutative}
\def\cpct/{compact}
\def\compl/{complex}
\def\comp/{complex}
\def\Comp/{Complex}
\def\conf/{conformal}
\def\conj/{conjugate}
\def\conn/{connect}
\def\cont/{continuous}
\def\conv/{converge} 
\def\convc/{convergence}
\def\convt/{convergent}
\def\convx/{convex}
\def\coord/{coordinate}
\def\lcoord/{local coordinate}
\def\Corr/{Corresponding}
\def\corr/{corresponding}
\def\corrd/{correspond}
\def\cov/{covariant}
\def\decomp/{decomposition}
\def\deco/{decompose}
\def\diff/{different} 
\def\Diff/{Different} 
\def\dimn/{dimension} 
\def\distr/{distribution} 
\def\div/{diverge} 
\def\dom/{domain}
\def\eg/{\hbox{\it e.g}\.}
\def\eigenf/{eigen\-\func/}
\def\eigensp/{eigen\-space}
\def\eigenv/{eigen\-value}
\def\eq/{equation}
\def\equa/{equation}
\def\de/{\diff/ial \equa/}
\def\do/{\diff/ial operator}
\def\ode/{ordinary \de/}
\def\pde/{partial \de/}
\def\pdo/{partial \diff/ial operator}
\def\psdo/{pseudo \diff/ial operator}
\def\fin/{finite}
\def\Ex/{\it Example.\ \normal}
\def\Exnr#1/{\it Example #1.\ \normal}
\def\foll/{follow}
\def\follg/{following}
\def\Follg/{Following}
\def\func/{function}
\def\Func/{Function}
\def\Fonc/{Fonc\-tion}
\def\fonc/{fonc\-tion}
\def\Funk/{Funk\-tion}
\def\funk/{Funk\-tion}
\def\gen/{general}
\def\har/{harmonic}
\def\Hint/{\it Hint. \normal}
\def\hist/{historic}
\def\histcl/{historical}
\def\hol/{holo\-morphic}
\def\homog/{ho\-mo\-ge\-ne\-ous}
\def\hyp/{hyper\-bolic}
\def\hyperg/{hyper\-geometric}
\def\ie/{\hbox{\it i.e.}}
\def\iff/{if and only if}
\def\Imm/{\operatorname{Im}}
\def\ineq/{inequality}
\def\infra/{{\it inf\-ra}}
\def\ultra/{{\it ult\-ra}}
\def\Inpart/{In particular}
\def\inpart/{in particular}
\def\instof/{instead of}
\def\interps/{interpolation space}
\def\interp/{interpolation}
\def\Interp/{Interpolation}
\def\interpr/{Interpretation}
\def\Intr/{Introduction}
\def\intv/{interval}
\def\inv/{invariant}
\def\invc/{invariance}
\def\Iowords/{In other words}
\def\iowords/{in other words}
\def\ipr/{inner product}
\def\irred/{irreducible}
\def\lb/{line bundle}
\def\lin/{linear}
\def\lhs/{left hand side}
\def\rhs/{right hand side}
\def\loc/{local}
\def\math/{mathematic}
\def\mathcn/{\math/ian}
\def\manif/{manifold}
\def\meas/{measure}
\def\measl/{measurable}
\def\mero/{mero\-morphic}
\def\mon/{monomial}
\def\monog/{monogenic}
\def\mult/{multiple}
\def\multy/{multiply}
\def\multn/{multiplication}
\def\nas/{necessary and sufficient}
\def\nbd/{neighborhood}
\def\neg/{negative}
\def\nondeg/{nondegenerate}
\def\Oohand/{On the other hand}
\def\oohand/{on the other hand}
\def\Oonhand/{On the one hand}
\def\oonhand/{on the one hand}
\def\oper/{operator}
\def\orth/{ortho\-gonal}
\def\orthon/{ortho\-normal}
\def\otoh/{on the other hand}
\def\quat/{quaternion}
\def\pp/{\hbox{a. e.}}
\def\psh/{plurisubharmonic}
\def\pol/{polynomial}
\def\pot/{potential}
\def\pos/{positive}
\def\princ/{principle}
\def\prob/{probability}
\def\proj/{projective}
\def\projn/{projection}
\def\Proof/{\it Proof:\normal}
\def\Rem/{\it Remark\normal}
\def\Ree/{\operatorname{Re}}
\def\Remnr#1/{\it Remark\ \normal #1. }
\def\rep/{representation}
\def\reps/{representations}
\def\meta/{metaplectic representation}
\def\repr/{reproducing}
\def\reprker/{reproducing kernel}
\def\resp/{respective} 
\def\resply/{respectively}
\def\restr/{restriction}
\def\sa/{self-adjoint}
\def\st/{such that}

\def\sol/{solution}
\def\ru/{space}
\def\sph/{spherical}
\def\ssp/{sub\ru/}
\def\sym/{symmetric}
\def\Sym/{Symmetric}
\def\symb/{symbol}
\def\symbc/{symbolic}
\def\symdom/{\sym/ domain}
\def\symp/{symplectic}
\def\Theor#1/{\fet Theorem #1.\ \normal}
\def\Lem#1/{\fet Lemma #1.\ \normal}
\def\Lemma/{\fet Lemma.\ \normal}
\def\topl/{topology}
\def\topll/{topological}
\def\transf/{transform}
\def\transl/{translation}
\def\transfn/{transformation}
\def\transv/{transvectant}
\def\trig/{trigonometric}
\def\tril/{trilinear}
\def\trilf/{trilinear form}
\def\uhp/{upper halfplane}
\def\uhs/{upper halfspace}
\def\vb/{vector bundle}
\def\vf/{vector field}
\def\vsp/{vector space}
\def\wrt/{with respect to}
\def\Wlog/{Without loss of generality}
\def\a{\alpha}

\def\sig{\sigma}

\def\Ab/{Abel}
\def\Ban/{Banach}
\def\Bansp/{\Ban/ space}
\def\Belt/{Bel\-tra\-mi}
\def\Berg/{Berg\-man}
\def\Bern/{Ber\-nou\-lli}
\def\Berz/{Berezin}
\def\Bess/{Bessel}
\def\Cart/{Car\-tan}
\def\Cay/{Cay\-ley}
\def\CG/{Clebsch-Gordan}
\def\Cl/{Clifford}
\def\CR/{Cauchy-Rie\-mann}
\def\Dir/{Dirichlet}
\def\Eucl/{Euclide}
\def\Eucln/{Euclidean}
\def\F/{Fourier}
\def\Hank/{Hankel}
\def\Hankf/{\Hank/ form}
\def\Herm/{Hermite}
\def\Hilb/{Hilbert}
\def\Hilbs/{Hilbert space}
\def\Hilbsp/{Hilbert space}
\def\HS/{Hilbert-Schmidt}
\def\Lag/{La\-grange}
\def\Lap/{La\-place}
\def\LapBelt/{\Lap/-\Belt/}
\def\Leb/{Lebesgue}
\def\Marc/{Mar\-cin\-kie\-wicz}
\def\Moeb/{Moebius}
\def\Moebt/{Moebius transformation}
\def\Moebtransfn/{Moebius transformation}
\def\Pla/{Plan\-che\-rel}
\def\Poin/{Poin\-car\'e}
\def\Riem/{Rie\-mann}
\def\Riemn/{\Riem/ian}
\def\psRiemn/{pseudo-\Riem/ian}
\def\Riems/{Rie\-mann surface}
\def\Schroe/{Schr\"odinger}
\def\Weier/{Weier\-strass}
\def\re{\operatorname{Re}}
\def\im{\operatorname{Im}}
\def\i/{\operatorname{\sqrt{-1}}}
\def\anal/{analytic}
\def\bsd/{bounded symmetric domain  }
\def\bdd/{bounded}
\def\calc/{calculation}\def\conj{conjugate}
\def\calci/{calculating}\def\eg{e.g.}
\def\conj/{conjugate}
\def\deco/{decomposition}
\def\eg/{e.g.}
\def\fct/{function}
\def\gp/{group}
\def\hw/{highest weight}
\def\hwv/{highest weight vector}
\def\hwvs/{highest weight vectors}
\def\lw/{lowest weight}
\def\lwv/{lowest weight vector}
\def\lwvs/{lowest weight vectors}
\def\hds/{holomorphic discrete series}
\def\iff/{if and only if}
\def\inv/{invariant}
\def\irrde/{irreducible decomposition}
\def\meas/{measure}
\def\transf/{transform}
\def\rep/{representation}
\def\resp/{respectively}
\def\inters/{intertwines}
\def\interg/{intertwining}
\def\meta/{metaplectic representation}
\def\qu/{quaternion}
\def\rep/{representation}
\def\symdom/{ symmetric domain}
\def\st/{such that}
\def\shd/{subhead}
\def\transf/{transform}
\def\wrt/{with respect to}

\def\Norm#1#2#3{\Vert#1\Vert^{#3}_{{#2}}}

\def\tr{\operatorname{tr}}

\begin{document}

\title[
  Principal series of Hermitian Lie groups 
induced from  Heisenberg parabolic subgroups
]
{
  Principal series of Hermitian Lie groups 
induced from  Heisenberg parabolic subgroups
}
\author{GenKai Zhang}
\address{Mathematical Sciences, Chalmers University of Technology and
Mathematical Sciences, G\"oteborg University, SE-412 96 G\"oteborg, Sweden}
\email{genkai@chalmers.se}

\thanks{ Research supported partially
 by the Swedish
Science Council (VR). 
}
\begin{abstract} Let $G$ be an irreducible Hermitian
  Lie group and $D=G/K$ its bounded symmetric domain
  in $\mathbb C^d$ of rank $r$.
Each $\gamma$ of the  Harish-Chandra
strongly orthogonal roots $\{\gamma_1, \cdots, \gamma_r\}$
defines
a Heisenberg parabolic subgroup
$P=MAN$ of $G$.
We study the 
principal series
  representations $\Ind_P^G(1\otimes e^\nu\otimes 1)$
  of $G$  induced   from $P$.
    These
  representations can be realized
  as the $L^2$-space on the minimal
  $K$-orbit $S=Ke=K/L$ of a root vector $e$ of $\gamma$ in $\mathbb C^d$, and $S$
  is a circle bundle over
  a compact Hermitian symmetric
  space $K/L_0$ of $K$ of rank one or two.
  We find the complementary series,  reduction points, and unitary
  subrepresentations   in this family of representations.
\end{abstract}
\maketitle

\baselineskip 1.40pc

\section{Introduction}

In the present paper we shall study composition series
and complementary series
for degenerate principal series representations for an irreducible Hermitian Lie group
$G$ induced from a Heisenberg parabolic subgroup.
Let $D=G/K\subset \mathbb C^d$ be the bounded
symmetric domain of $G$ in its Harish-Chandra
realizaation in $\mathbb C^d=\fp^+$.
Any choice of a Harish-Chandra
strongly orthogonal root determines a one-dimensional
split subgroup $A=\mathbb R^+$ of $G$ and 
 a Heisenberg parabolic subgroup
$P=MAN$ of $G$. 
We study 
the induced representation
$I(\nu)=
\Ind_P^G(1\otimes e^\nu\otimes 1)$ of $G$ from $P$.  We shall find its
complementary
series, the reduction points, explicit realization
of certain finite dimensional representations,
and unitarizable subrepresentations.

The representation 
 $I(\nu)$ can be realized
on the $L^2$-space on a homogeneous space $K/L$ of $K$, $L=M\cap K$,
and $K/L$ is an orbit of $K$ in $\mathbb C^d$, with the Harish-Chandra
root vector as a base point. It is a circle bundle $K/L\to K/L_0$ over 
 its projectivization $K/L_0$.
The space $K/L_0$ itself
is a {\it compact Hermitian symmetric space}. We can find the
irreducible decomposition of $L^2(K/L)$
by using the Cartan-Helgason
theorem for line bundles over $K/L_0$. We study
then the Lie algebra action of $\fg$ on 
 $L^2(K/L)$.

The induced representations $Ind_P^G(\nu)$
from Heisenberg parabolic subgroup $P$
can be viewed as the counterpart of
the representations $Ind_Q^G(\nu)$
from the Siegel parabolic subgroup $Q$, both $P$ and $Q$ being maximal
parabolic subgroups.
The nilpotent part  in $P$
is a Heisenberg nilpotent group  whereas
it is abelian in $Q$. 
The representations 
$Ind_Q^G(\nu)$ and their analogues 
are well studied and 
are of considerable interests 
as they are closely related to the holomorphic discrete 
series \cite{OZ-duke}. 
The general case of maximal parabolic subgroups
with abelian nilradical can be put in the
setup of Koecher's construction of Lie algebras
and the corresponding induced representations
have also been intensively studied; see
e.g. \cite{Sahi-crelle, Z-matann}.
The analysis is done in \cite{Sahi-crelle} by using
eigenvalues of intertwining differential operators
and the tensor product structure of induced representations,
and in \cite{Z-matann}
by computing differentiations
and recursions of spherical polynomials along a torus in the symmetric space
$K/L$.
However in our case the
nilpotent group $N$ is non-abelian and $K/L$ is not symmetric, there is no known
construction of intertwining  operators, so the
method in \cite{Sahi-crelle}  seems not possible,
and we shall adapt the method in  \cite{Z-matann}.
Now the manifold $K/L_0$ is a complex manifold
and the differentiation  involves vector fields in $T^{(1, 0)}(K/L_0)$,
we have to develop further the technique in loc.~cit.~.
We consider the differentiation along products of projective
sphere $SU(2)/U(1)$ in $K/L_0$ and find the differentiation formulas
by considering classical spherical polynomials of $SU(2)$.
By using Weyl group symmetry of spherical functions we
find then the required formulas for  the Lie algebra actions.

The study of induced
representations from Heisenberg parabolic
subgroups can be put in a rather general context.
It has been proved by Howe \cite{Howe} 
that in most  semisimple Lie groups $G$ 
there are Heisenberg parabolic subgroups, and these groups  
have been all classified.  
The introduction of them is closed related
to Howe's notion of the rank and minimal representations.
Recently Frahm \cite{Frahm}
has studied the intertwining differential operators
for the induced representations for  Heisenberg parabolic subgroups
in some non-Hermitian Lie groups.
Also spherical duals of classical groups have
been studied extensively; see \cite{Barbasch, KS, KS-2}.
When $G=SU(p, q)$ the induced representations
can be realized on a homogeneous cone
in $\mathbb C^{p+q}$ and they have been studied
by Howe-Tan \cite{HT} in greater details.

We mention that
the space $K/L_0$ has an interesting geometry, it is
 the variety of minimal rational tangents
in a fixed tangent space $V=T_0^{(1, 0)}(D)$ of the symmetric space $D$; more precisely 
it is the projectivization $\mathbb P(V)$ of all tangent vectors $v\in V$
with maximal holomorphic sectional curvature. It is also
projectivization of the  space $K/L$ of minimal tripotents
and plays important role in complex differential geometry; see \cite{Hwang-Mok, Mok-geo-st}.

The paper is organized as follows.
In Section 1 we introduce the parabolic
subalgebra $\fp=\fm +\fa +\fn=\fm +\mathbb R\xi +\fn$ and the principal series
$(I(\nu), \pi_\nu)$.
In Section 2 we find the irreducible decomposition
for $I(\nu)=L^2(G/P)=L^2(K/L)$ under $K$.
The action of $\pi_\nu(\xi)$ on 
$I(\nu)$  is done in Section 3 and it is one of our main results.
As consequences we find the complementary series, certain unitarizable subrepresentations,
and also realization of certain finite dimensional
representations of $\fg^{\mathbb C}$ in the induced representations
in Section 4. In Section 5 we treat the case when
$\fg=\fsu(d, 1), \fsp(r, \mathbb R)$
where $K/L$ is the sphere and $K/L_0$
is the complex projective space. We give a simpler proof
of the classical result of Johnson-Wallach \cite{JW}.
In Appendix \ref{app-a}
 we compute certain
recurrence formulas for spherical polynomials
over the complex projective line $\mathbb P^1$;
we need only the leading coefficients in the formulas
which can be easily proved by other methods, but we present
the complete formulas as they might be of independent interests
in special functions.  The complete list of the spaces $G/K$ and $K/L_0$
are given in Appendix \ref{app-b}.

I would like to thank Jan Frahm for  discussions
on Heisenberg parabolic subgroups. I'm also grateful
to Dan Barbarsch for
explaining his work
on unitary spherical duals
of classical groups and Harald Upmeier for explanation
of Jordan quadrangles.

\subsection*{Notation}

 Hermitian
 symmetric spaces
 have rather rich structure so is
 the notation. For the convenience of the reader we list
the main  symbols used in our paper.

\begin{enumerate}
  \item
  $D=G/K$, 
  bounded symmetric domain of $D$
  of rank $r$
  realized in the Jordan triple $V=\mathbb C^d$
  with Jordan triple product $\{x, \bar y, z\} =D(x,  y) z$
  and Jordan characteristic $(a, b)$ (or root multiplicities);
  $    d=r+ \frac a2 r(r-1) +rb $, the  dimension.

\item     $\fg=\fk+\fp$,
  Cartan decomposition, 
  $\fp=\{\xi(v) = v+Q(z) \bar v; v\in V\}$
  as  holomorphic vector fields on $D$.

\item   $\fg^{\mathbb C} =
  \fp^{-}
+  \fk^{\mathbb C }
  +    \fp^{+} 
  $,
 Harish-Chandra decomposition
  respect to the center element $Z$
  of $\mathfrak k$, $\Ad (Z)|_{\fp^{\pm}} =\pm i$;
  $\fp^+ = V$.

\item $e=e_1$, a fixed minimal tripotent, $V=V_1+V_1 +V_0=
  V_1(e)+V_1(e) +V_0(e),
  $
  the Peirce decomposition with respect to $e$,
  and
  $\{e, v_1, w, v_2\}$
  a Jordan quadrangle.

\item $\{e, \bar e, D(e,  e)\}$,
  standard $\mathfrak{sl}(2)$-triple; $H_0=iD(e,  e)$.

\item $\xi=\xi(e)=e+\bar e\in \fp$, $\fa=\mathbb R\xi$,
  $\fg=\fn_{-2}+
  \fn_{-1}+ \fm +\fa
  + \fn_1 +\fn_2
  $,  the root space
  decomposition of $\fg$ with respect to $\fa$,
 $\fn=
  \fn_1 +\fn_2$,  a Heisenberg Lie algebra;
$\rho_{\fg}=1+(r-1) a +b$, half sum of positive roots.

\item $M$, $A$, $N$, the corresponding subgroups
  with Lie algebra $M$, $A$, $N$, 
  $L=M\cap K=\{k\in K; k e=e\}\subseteq K$
  with Lie algebra $\fl$.

\item  $\fk_1=[\fk, \fk]$, the semisimple part 
of $\fk$, and $K_1\subset K$ the corresponding Lie group.

\item   $S=K/L=G/P$, the  manifold of
  rank one tripotents.
\item
$ S_1  =\mathbb P(S)=K/L_0=K_1/L_1 $, the projective
  space of $S$, also called
  the variety of minimal rational tangents;
  $L_0=\{k\in K; k e=\chi(k) e, \chi(k)\in U(1)\}$
  the subgroup of elements of $K$ fixing the line $\mathbb Ce$
  with Lie algebra   $\fl_0=\{X\in \fk; Xe=\chi(X) e, c\in i\mathbb R\}$.
  $ S_1=K/L_0=K_1/L_1$,  semisimple  Hermitian symmetric space
  of rank two of dimension $d_1=\dim_{\mathbb C} V_1=(r-1) a +b$,
  $  (a_1, b_1)$, the Jordan
  characteristic of $S_1$;
  
\item    $\chi_l(k)=\chi(k)^l$, character on $L_0$, $l\in \mathbb Z$.
  
\item
  $\fk=\fq +\fl_0$, $\fk_1=\fq +\fl_1$
  Cartan decomposition
  for the symmetric space
 $ S_1=K/L_0=K_1/L_1$;
  $\fq= \{ D(v,   e) -D(e,   v), v\in V_1\}$;
  $\fl_1=\fl_0\cap \fk_1$, the semisimple component of
  $\fk_0$.
\item $\fk^{\mathbb C}
= \fq^{-}  +
  \fl_0^{\mathbb C} + 
   \fq^{+}
   $,
   $\fk_1^{\mathbb C}
= \fq^{-}  + \fl_1^{\mathbb C} + 
   \fq^{+}
   $,    the Harish-Chandra decomposition of $\fk^{\mathbb C}, \fk_1^{\mathbb C}$
   for the Hermitian symmetric space $K/L_0=K_1/L_1$;  
 $\fq^+= \{ D(v,  e);  v\in V_1\}$;
$\fq^-= \{ D(e,   v);  v\in V_1\}$;
the Jordan triple defined  on $\fq^+= \{ D(v,  e);  v\in V_1\}$
is via the Lie bracket,  $$
[[D(v, e), D(e, w)], D(v, e)] = D(D(v, w)u,
e), $$ and is isomorphic to the Jordan triple $V_1\subset V$;
$\fq^+$
is the
 the holomorphic tangent space  
 of $S_1=K/L_0=K_1/L_1$ at $L_0\in K/L_0$.

  \item $\fk_1^\ast=\fl_1 + i\fq$, the non-compact
   dual of $\fk_1=\fl_1 + \fq$;
   $\fh_{i\fq}
   \subset  i\fq$, Cartan subspace.

\item $\rho:=  \rho_{\fk^\ast_1}:=\rho_1\alpha_1 
+ \rho_2 \alpha_2=(\rho_1, \rho_2), \quad \rho_1=1+a_1+b_1, \, \rho_2=1+b_1,$
half sum of positive restricted roots of $\fk^*_1$
with respect to     $\fh_{i\fq}$.
  
 \item  $\Phi_{\lambda, l}$, Harish-Chandra spherical
   function for the symmetric
   pair  $(\fk_1^\ast, \fl_1)$ 
   with one-dimensional character  $\chi_l$, 
$C(\lambda, l)=C(\lambda, -l)$, Harish-Chandra $C$-function 
  for the symmetric space $S_1=K/L_0=K_1/L_1$ with  character 
    $\chi_l$; $\phi_{\mu, l}$ spherical polynomial.
  
\end{enumerate}  

\section{Preliminaries}
We shall use the Jordan triple description of Hermitian symmetric
spaces; see
\cite{Loos-bsd,  Upmeier}.

\subsection{Hermitian symmetric space $D=G/K$
}
\label{1.1}
Let $D$ be an irreducible bounded symmetric domain of rank $r$ in  $V=\mathbb
C^d$.
Let $G$ be the group of
bi-holomorphic automorphisms
of $D$,  and $K=\{k\in G; k0=0\}$
the maximal compact subgroup of $G$, so
that $D=G/K$.
The space $V$ has the structure
of an irreducible Jordan triple with triple product
$\{x, \overline y, z\}=D(x,  y)z$ with the corresponding
$\text{End}(\bar V,  V)$-valued quadratic
  form $Q(x)$, $Q(x)\overline y=\frac 12 D(x,   y)x$.
  Note that
  in \cite{Loos-bsd} $D(x,  y)$ 
  is written as $D(x,  \overline y)$,
  and to  ease notation
  we write it just as   $D(x, y)$
 so it is {\it conjugate linear} in $y$.
Let   $(a, b)$ be the Jordan characteristic 
 of $V$,  and $b=0$ when $D$ is of tube domain. 
 The   dimension  $  d=r+ \frac a2 r(r-1) +rb $.

   Let $\fg=\fk +\fp$
  be the Cartan decomposition
  of $\fg$. Realized as holomorphic vector fields
  on $D$, i.e., as $V$-valued functions on $D$, the space
  $\fp$ is
  \begin{equation}
    \label{cartan-p}
  \fp=\{\xi(v)=v-Q(z) \overline v; v\in V\}.    
  \end{equation}
  The adjoint action $v\mapsto \ad(k)v$  of $K$
  as well as $\fk$
on $\fp$ coincides with its defining 
action on $D$ and will be written just as $kv=\ad k\, v$,
$Xv=\Ad(X) v$, $k\in K$, $X\in \fk$ when no confusion would arise.

Denote $Z\in \fk$  the central element defining
  the complex structure of $\fp$ and
$\fg^{\mathbb C}=\fp^+ +
\fk^{\mathbb C} +\fp^-$ be the Harish-Chandra decomposition,
$Z|_{\fp^{\pm}} =\pm i$. The spaces $\fp^{+}$ is identified  with $V$
via the identification $V\ni v =\frac{1}2 (\xi_{v} -i \xi_{iv}) \in \fp^+$
and $\overline V=\{\overline v; v\in V\}$ with $\{-Q(z)v\}=\fp^-$.
The Lie algebra $\fk=\fk_1\oplus \mathbb R(iZ)$,
where $\fk_1=[\fk, \fk]$ is the semisimple part of
$\fk$ with trivial center. Let $K_1\subset K$ be the corresponding
semisimple subgroup of $K$
with Lie algebra $\fk_1$.

We fix  the Euclidean inner product on $V$ so that 
a minimal tripotent has norm $1$, and
fix the corresponding normalization
of the Killing form on $\fg^{\mathbb C}$.
All orthogonality in Lie algebra
 $\fg^{\mathbb C}$ below 
is with respect to the Killing form.

\subsection{Maximal parabolic subgroup $P=MAN$ of G
  and induced representation  $Ind_P^G(  \nu)  $ }
We fix in the rest of the paper a minimal tripotent 
 $e=e_1$ and denote
 \begin{equation}
   \label{xi-H0}
   \xi=\xi_e,\quad H_0= iD(e,  e),
\quad \fa=\mathbb R\xi\subset 
\fp. \end{equation}
A Harish-Chandra strongly orthogonal
root $\gamma_1$ for $\mathfrak g^{\mathbb C}$
can be chosen so that its co-root is $D(e, e)$, $\gamma_1(D(e, e))=2$.
We shall only need $\gamma_1$ below.


The Peirce decomposition of $V=\mathbb C^d$ with respect to
the  tripotent $e$ is
\begin{equation}
  \label{Peirce}
V=V_2 + V_1 + V_0, \, V_j=V_j(e):= \{v\in V; D(e,  e)v=j v\}, \, j=0, 1, 2,
\end{equation}
Furthermore  $V_2=\mathbb C e$
is one dimensional,
$V_1$ is of dimension  $d_1=\dim_{\mathbb C} V_1=(r-1) a +b$,
and $V_0$ is a Jordan triple of rank $r-1$
and dimension $ 1 +\frac 12 a(r-1)(r-2) + (r-1)b$.
The Jordan rank of $V_1$  is
\begin{equation}
  \label{rank-V-1}
  \rank V_1=\begin{cases} 2, & \fg\ne \fsu(d, 1), \fsp(m,  \mathbb
    R)\\
    1,  &\fg= \fsu(d, 1), \fsp(m,  \mathbb R).
  \end{cases} 
\end{equation}
Certain computations have to be done
depending on  the different cases.

We shall need the description for the root spaces of
$\fg$ under $\fa=\mathbb R\xi$.

\begin{lemm+}
  \label{root-1}
  The root space decomposition of
  $\fg$ under $\fa=\mathbb R\xi$ is
\begin{equation}
  \label{nman}
  \fg =\fn_{-} +
  (\fm +\fa)+
    \fn, \quad 
\fn=\fn_{1}+\fn_{2}, \quad \fn_-=\fn_{-1}+\fn_{-2},
\end{equation}
where $\fn_{\pm 2}$, $\fn_{\pm 1}$, 
and $\fm +\fa$ are the root spaces of $\xi$ with 
roots $\pm 2, \pm 1, 0$, respectively. The subspaces are given by
\begin{equation}
  \label{n0}
  \fm =\fl \oplus \{\xi_v; v\in V_0\}, \quad 
  \fl =\fm \cap \fk =\{X\in \fk; Xe=0\}, 
\end{equation}
\begin{equation}
  \label{n1}
  \fn_{1}
  =\{\xi_{v} +(D(e,  v) -
D(v,  e )); 
v\in V_1\}, 
\end{equation}
and 
\begin{equation}
  \label{n2}
\fn_{2}=\mathbb R(\xi_{ie} - H_0). 
\end{equation}
The half sum of positive roots 
is \begin{equation}
\rho_{\fg}:=1+(r-1) a +b=1+\dim_{\mathbb C} V_1.\end{equation}
\end{lemm+}  
\begin{proof}
  The root spaces are described in  \cite[Lemma 9.14]{Loos-bsd}.
  The  root space $\fn_{1}$ has real dimension
  $2  \dim_{\mathbb C} V_1
  $ and thus
  $\rho_{\fg} = 1 +  \dim_{\mathbb C} V_1 = 1+ (r-1)a+ b$.
\end{proof}

The nilpotent algebra $\fn$ is a Heisenberg Lie
algebra. Their appearance in general semisimple Lie algebras
has been classified in \cite{Howe}.

Let $M, A, N$, $L=M\cap K=\{k\in K; k e=e\}$,  be the corresponding Lie subgroups
and $P=MAN$ the parabolic group of $G$.

A linear functional $\nu$ on $\fa^{\mathbb C}$ will be identified as
$\nu\in
\mathbb C, \nu=\nu( \xi)$.  The main object of this paper is 
the induced representation 
\begin{equation}
 \label{def-ind}
 I(\nu)
 =Ind_P^G(\nu):=Ind_{P}^{G}
(1\otimes e^{\nu}\otimes 1) 
\end{equation}
defined as the space of measurable functions 
on $G$ such that 
$$
f(gme^{t\xi}n)=  e^{-t\nu} f(g), m\in M, n\in N, t\in 
\mathbb R, 
$$
and $f|_{K}\in L^2(K)$. Any $f
\in
I(\nu)$ as function on  $K$ is right
$L$-invariant, and thus $f\in L^2(K/L)$. The corresponding
$(\fg^{\mathbb C}, K)$-representation,   $\fg^{\mathbb C}$ being
the complexification
$\fg$,  will also be denoted by $I(\nu)$.
If $\nu = \rho_{\fg} + i\lambda, \lambda\in \mathbb R, \lambda\ne 0$, then
$I(\nu)$ is a unitary irreducible
representation of $G$ on $L^2(K/L)$; see e.g. \cite{Kn}.

\section{Decomposition of 
  $(L^2(K/L), K)$ and 
  spherical polynomials}

We assume throughout Sections 2 and 3 that the Hermitian 
symmetric symmetric domain  $D=G/K$
is irreducible of rank $r\ge 2$ 
and  is not the Siegel domain 
$Sp(r,  \mathbb R)/U(n)$; this case
and the rank one domain
$SU(d, 1)/U(d)$
will be treated in  Section 5.

\subsection{The homogeneous space 
  $S=G/P=K/L$ as circle bundle over 
compact   Hermitian symmetric space $S_1=K/L_0=K_1/L_1$}

The homogeneous space $S:=G/P=K/L$ of $K$
can be realized as the manifold 
of rank one 
tripotents in $V$, $S=K e\subset V$. We shall
fix this realization in the rest of the paper and  use
the global coordinates on $V$ for $S$ when needed.

To find the decomposition
of $L^2(K/L)$ under $K$ we consider the projectivization
$$
S\to S_1=\mathbb P(S)
=\{[v]:=\mathbb C v\in \mathbb P(V); v\in S\}, \quad
v\mapsto [v]$$
of  $S$ in the projective space 
$\mathbb P(V) $
of $V$. 
Then $S_1=K/L_0$,
where $L_0=\{k\in K; k e\in \mathbb Ce\}$,
and $L\subset L_0$ is a normal subgroup with $L_0/L$
being the circle group $U(1)$.
The natural map $S=K/L\to  S_1=\mathbb P(S)=K/L_0$ defines
a fibration
\begin{equation}
\label{circle-bdl}  
S=K/L\to S_1=K/L_0, 
\end{equation}
of $\mathbb P(S)$
with fiber  the circle $U(1)=L_0/L$.

The space  $S_1=\mathbb P(S)$  is a compact Hermitian symmetric space of (non-simple) $K$.
The complete list of $(G, K, L_0)$, $D=G/K$, $S=G/P=K/L\to S_1=K/L_0$
is given in Tables  \ref{tab:1}-\ref{tab:2}. 
The space $S_1=K/L_0$ is also the  variety of minimal rational tangents (VMRT) \cite{Hwang-Mok} of 
the symmetric space $D$.

The action of $L_0$ on $e$
defines a character of $\chi: L_0\to U(1)$, namely 
\begin{equation}
  \label{def-chi}
  L_0=\{k\in K; k e\in \mathbb Ce\}=  \{k\in K; k e=\chi(k)e\}.
\end{equation}

The spaces $S$ and $S_1$
are also homogeneous space of the semisimple
part $K_1\subset K$. Let $L_1=\{k\in K_1; k e =\chi(k)e\}\subset L_0$,
then $S_1=K/L_0=K_1/L_1$, $S=K_1/L\cap L_1$.
The fibration (\ref{circle-bdl}) above becomes
a circle bundle $S=K/L\to S_1=K/L_0   =K_1/L_1$
for   $K_1$-homogeneous spaces with
the fiber  $U(1)=L_1/L\cap L_1$.

The element $\exp(\pi  H_0), H_0= iD(e, e),$
defines 
a Cartan involution $\Ad\exp(\pi  H_0)$ on $\fk$
with the corresponding Cartan decomposition 
$$
\fk 
=\fl_0 +\fq.$$
The space $\fq$, similar to (\ref{cartan-p}), 
is represented  by 
\begin{equation}
  \label{cartan-q}
  \fq=\{D(v,   e)-D(e,  v); 
v\in V_1 
\}.
\end{equation}

We fix now a
complex structure on $\fq$
and the corresponding Harish-Chandra decomposition
of $\fk^{\mathbb C}$.

\begin{lemm+}\label{cplx-q}
  Define the $K$-invariant
  complex structure on $S_1=K/L_0$ by
  the element $iD(e, e)\in \fl$, i.e. by
  $\frac 12 \Ad (iD(e, e))\in \fl$ on $\fq\fq=T_{[e]}(S_1)$,
  at the base point $[e]=\mathbb Ce S_1$.
  In terms of elements in $\fl_1'$
in the semi-simple part $\fk_1^{\mathbb C}$ 
of $\fk_1^{\mathbb C}$, it is the same
as
\begin{equation}
  \label{semi-simple}
-\frac 12 \Ad (H_0)|_{\fq}
=\frac 12 \Ad (Z')|_{\fq}, \quad 
Z'=i\frac{p}{d} Z- H_0\in \fl_1\subset \fk_1.  
\end{equation}
The Harish-Chandra decomposition of $
\fk^{\mathbb C}$ is 
\begin{equation}
  \label{hc-k}
  \fk^{\mathbb C}
=\fq^{-} +
\fl_0^{\mathbb C} +\fq^{+}
\end{equation}
with 
\begin{equation}
  \label{q+-}
T_{[e]}^{(1, 0)}(S_1)=\fq^{+}
=\{D(v, e); 
v\in V_1 
\}, \, 
T_{[e]}^{(0, 1)}=\fq^{-}
=\{-D(e, v); 
v\in V_1 
\}.  
\end{equation}
\end{lemm+}
\begin{proof} The elements $cZ, c\in \mathbb R$ is in the center
  of $\fk$ so that $\Ad (X)=
  \Ad (cZ+X)$
  on $\fk$  for any $X\in \fk$.
  The semisimple component of $X$
  in $\fk_1$  is obtained
  as $X-   \frac 1d \tr \Ad (X)|_{\fp^+}$, $d=\dim V$.
  Thus we have   (\ref{semi-simple}). The Harish-Chandra
  decomposition is obtained by the commutator
  formula in Jordan  triples   \cite{Loos-bsd},
  $$
  [D(u, v),   D(x, y)]
  = D(D(u, v)x, y) - D(x, D(v, u)y).
  $$
\end{proof}  

Here  the complex structure on $\fq$ is chosen
so that $$
v\in V_1\subset V=\fp^+\to D(v, e)\in \fq^+$$
is complex linear so that
the complex structures in
$\fp$ and $\fq$ match in this sense.
The Lie algebra structure in $\fk^{\mathbb C}$
defines $\fq^{\pm}$
as a Jordan triple, 
and it is isomorphic to  $V_1$. To avoid confusion we shall keep the notation 
$\fq^{+}$.

\subsection{Cartan subalgebra $\fh_{i\fq}  \subset i\fq$
  and the restricted 
  root system for the non-compact 
  symmetric pair $(\fk^\ast, \fl_0)  =( \fl_0 + i\fq,  \fl_0)$}
We construct now split Cartan subalgebra 
in $i\fq$ for the
symmetric pair
$(\fk^\ast, \fl_0)  =( \fl_0 + i\fq,  \fl_0)$
and its semisimple part
$(\fk^\ast_1, \fl_1)  =( \fl_1 + i\fq,  \fl_1)$,
and find the corresponding root system. They will be used
in the decomposition of $L^2(K/L)$
and computation
of Harish-Chandra $c$-function.
We need the notion of a Jordan {\it quadrangle}
l \cite[p.~12, p.16]{N}.

An ordered quadruple $(u_0, u_1, u_2, u_3)$
of minimal tripotents is called 
Jordan quadrangle if the following three conditions
are satisfied,
 for all $i$ modulo $4$,
\begin{enumerate}
  \item $u_i$ and $u_{i+1}$
are in each other's Peirce $V_1$-space, 
$u_i\in V_1(u_{i+1})$, 
$u_{i+1}\in V_1(u_{i})$;
\item
   $u_i$ and $ u_{i+2}$
   are orthogonal as tripotents;
\item
$D(u_i,  u_{i+1})u_{i+2} = u_{i+3}$.
\end{enumerate}
Recall that we have assumed in this section that 
the domain $D\ne III_n=Sp(n, \mathbb R)/U(n),
D\ne I_{n, 1}=SU(n, 1)/U(n)$  (the Type IV domain $IV_{3}= III_2$ is
also excluded).
Then starting with the fixed minimal tripotent
$e$ there are minimal tripotents $v_1, w, v_2$
such that 
$(u_0, u_1, u_2, u_3)=(e, v_1, w, v_2)$ is a Jordan quadrangle.
This is implicitly 
in  \cite{N} where orthogonal bases (called grids)
are constructed for Jordan triple systems,
and we provide brief arguments. The Jordan triple
system $V$ is of rank $r\ge 2$ so there exists
a Jordan algebra  $V'$ as a sub-triple of $V$
of rank two with $e_1+e_2$ as identity element,
where $e_1=e, e_2$ are the Harish-Chandra strongly
orthogonal root vectors;
for  $$D=I_{r, r+b}, \,  II_{2r}, \, II_{2r+1}, \, IV_{n} (n>4), \, V, \,
VI$$
the corresponding $D'$ is
$$
D'=I_{2, 2},  \, II_{4},  \, II_{4}, \, IV_{4}, \, IV_{8},  \, IV_{10}.$$
In all cases the  Jordan algebra $M_{2, 2}$ 
of square $2\times 2$-matrices, $D=I_{2,2}=IV_4$, forms a Jordan sub-triple, since
$$I_{2, 2}=IV_4\subset  II_{4}=IV_{6}\subset  IV_{8}\subset IV_{10}$$
in the sense of  Jordan sub-triples.
The following standard matrices
$$
E_{11}, E_{12}, E_{22}, E_{21} 
$$
form a  Jordan quadrangle in
$M_{2, 2}$ and in $V$.
 
Fix in the rest of the paper the Jordan quadrangle $\{e, v_1, w,
v_2\}$.
We have 
\begin{equation}
H_0 v_1=iD(e,  e)v_1=iv_1, \, H_0 v_2=iD(e,  e)v_2, D(v_1, 
v_1)e=D(v_2,  v_2)e=e,
\end{equation}
and 
  \begin{equation}
  \label{quadruple}
  D(e,   w) =D(v_1,  v_2)=0, \, D(v_1,  e)v_1=
  D(v_2,  e)v_2=0, 
  D(e,  {v_1})w =v_2, \, 
D(v_1,   e) v_2= w, 
\end{equation}
which we shall use below. 
See \cite[p.~12, p.16]{N} for further details.


The above construction results in the following two commuting copies of $sl(2, \mathbb C)$-triples in $\fk_1^{\mathbb C}$, 
\begin{equation}
  \label{H-12}
  E_j^+=D(v_j, e)\in \fq^+,  \, E_j^-=D(e, v_j)\in \fq^-, \,
  \, H_j= D(v_j,  v_j)-D(e,  e)\in \fk^{\mathbb C},\,
  \ E_j:= E_j^+- E_j^- \in \fq,
\end{equation}
with  the canonical relation
$$[H_j, E_j^{\pm}]= 2 E_j^{\pm},\, 
[E_j^+, E_j^-]= H_j, \, j=1, 2. 
$$
Moreover $\mathbb R (iE_1) +\mathbb R(iE_2)\subset i\fq $ is maximal abelian.
\begin{defi+}
  Let
\begin{equation}
  \label{cartan-h}
{
      \fh_{\fq}=
    \mathbb R E_1 +\mathbb R E_2 \subset \fq, \quad
    \fh_{i\fq}=
    \mathbb R (iE_1) +\mathbb R(iE_2)\subset i\fq 
  }, \quad
  \fh_{\fq}^{\mathbb C}=
{\mathbb C} E_1 +\mathbb C E_2.
\end{equation}
Extend the abelian subalgebra
$\mathbb CE_1 +\mathbb CE_2
$ of $\fk_1^{\mathbb C}$ 
to a Cartan subalgebra 
$$
\fh^{\mathbb C}_1:=  (\mathbb CE_1 +\mathbb CE_2) \oplus \mathfrak h_+, \quad 
\mathfrak h_+\subset \fl_1^{\mathbb c}\subset \fk_1^{\mathbb C}, 
$$
of $\fk_1^{\mathbb C}$,
so that 
$$
\fh^{\mathbb C}:=
 (\mathbb C Z+\mathbb CE_1 +\mathbb CE_2) \oplus \mathfrak h_+
$$
is a Cartan subalgebra of $\fk^{\mathbb C}$ and $\fg^{\mathbb C}$.
Define $\{\alpha_0, \alpha_1, \alpha_2\}$ to be the dual 
basis vectors of $\{iZ, iE_1 , iE_2\}$
that are vanishing  on $\mathfrak h_+$. 
\end{defi+}

Now
\begin{equation}
  \label{cartan-f-all}
(\mathbb C Z+\mathbb CH_1 +\mathbb CH_2) \oplus \mathfrak h_+\subset 
\fl_0^{\mathbb C}\subset \fk^{\mathbb C}\subset \fg^{\mathbb C}
\end{equation}
is  a Cartan subalgebra of three algebras
$\fl_0^{\mathbb C}, \fk^{\mathbb C}$ and $\fg^{\mathbb C}$.

We shall need the Cartan-Helgason theorem in \cite{Sch}
for line bundles over $K/L_0$ defined
by the  characters $\chi_l=\chi^l$; the character in \cite{Sch} 
is define using Cayley transform and Cartan subalgebras instead 
the geometric definition here. The relevant Cayley transform in our setup is
$$
c=c_{\fk}=\ad \(\exp(-\frac{\pi i}{4}(E_1^++E_2^+ + E_1^-+ E_2^-)\).
$$

\begin{lemm+}
  \label{k-roots}
  \begin{enumerate}
  \item  
  The subspace
  $\fh_{i\fq}$
  is maximal  abelian in $i\fq$.
If  $ \fg\ne \fsu(r+b, r)$, $r>1$
then  $K/L_0$ is an irreducible symmetric space 
and the restricted root system 
  for the non-compact  dual 
  $\fk^*=\fl_0 + i q$ 
with respect to $\fh_{i\fq}$ is 
\begin{equation}
  \label{root-sys-q}
  R(\fk^*, \fh_{i \fq})=  R(\fk^*_1, \fh_{i \fq}):
  =  \{\pm 2\alpha_1,\, \pm 2\alpha_2\}\cup 
  \{\pm \alpha_1\pm \alpha_2\}\cup 
 \{ \pm \alpha_1,\, \pm \alpha_2\}  
\end{equation}
with root multiplicities $(1, a_1, 2b_1)$ for the three  subsets of 
roots, $a_1, b_1$ being given in Tables \ref{tab:1}-\ref{tab:2}. 
The half-sum of the positive roots with respect to the ordering 
$\alpha_1>\alpha_2>0$ 
is 
\begin{equation}
  \label{rho-g}
\rho:=  \rho_{\fk^\ast_1}=\rho_1\alpha_1 
+ \rho_2 \alpha_2=(\rho_1, \rho_2), \quad \rho_1=1+a_1+b_1, \, \rho_2=1+b_1. 
\end{equation}
The two linear functionals $\rho_{\fk^\ast_1} $ and $\rho_{\fg}$ are related 
by
$$\rho_{\fg}=1+\rho_1 +\rho_2.$$

If  $ \fg= \fsu(r+b, r)$, $r>1$
then  $K/L_0=K_1/L_1$ is  reducible, $\fk_1=\fsu(r+b) +\fsu(r)$,
$\fk_1^\ast=\fsu(1, r+b-1) +\fsu(1, r-1)$, and the restricted root system 
is
\begin{equation}
  \label{root-sys-q-su}
R(\fk^*, \fh_{i q})=  R(\fk^*_1, \fh_{i q}):
  =  \{\pm 2\alpha_1,  \pm \alpha_1\} \cup
  \{  \pm 2\alpha_2,   \pm \alpha_2\}  
\end{equation}
with root multiplicities $(1, r+b-1), (1, r-1)$ respectively.
The corresponding $\rho$
is 
\begin{equation}
  \label{rho-k}
  \rho_{\fk^\ast_1}
  =\rho_1\alpha_1 
+ \rho_2 \alpha_2=(\rho_1, \rho_2), \quad \rho_1=r+b-1, \, \rho_2=r-1
\end{equation}
The relation  $\rho_{\fg}=1+\rho_1 +\rho_2$ holds also in this case.

\item
  The Cayley transform  $c$ exchanges two Cartan subalgebras
  $$
  c:
(\mathbb C Z+\mathbb CH_1 +\mathbb CH_2) \oplus \mathfrak h_+
  \to \fh^{\mathbb C},
  $$
$$c: H_j\to   iE_j$$
and the  pullbacks 
$$
\text{$c^\ast(2\alpha_j), j=1, 2,$  are the Harish-Chandra  orthogonal roots 
for  $(\fk_1^\ast, \fl_1)$.}
$$
\end{enumerate}\end{lemm+}

\begin{proof} We have assumed that
  $D$ is not $SU(d, 1)/U(d)$ nor $Sp(n, \mathbb R)/U(n)$,
so that $K/L_0=K_1/L_1$ is of rank two; see the Tables \ref{tab:1}-\ref{tab:2}. 
The rest  is a consequence 
of general results on 
root systems of non-compact Hermitian 
Lie algebra applied to the 
non-compact Helgason dual $\fk^*=\fl_0 + i q$
of the compact $\fk=\fl_0 + i q$. It follows 
also from  \cite[Lemma 
9.14]{Loos-bsd}, the Lie algebra
$\fg=\fp +\fk$ there being replaced by
our $\fk^\ast_1=i\fq +\fl_1$.
The abelian subspace in $i\fq$
is obtained from
the Harish-Chandra root vectors corresponding
to a frame of  minimal tripotents in $\fq^+$, 
and in the present case the frame in $\fq^+$
is $\{D(v_1, e), D(v_2, e)\}$.
Finally  the dimension of a general
Jordan triple of characteristics $(r, a, b)$
is $r +\frac 12 r(r-1)a+rb$. We have
 by Lemma   \ref{root-1}, 
$\rho_{\fg}=1+\dim_{\mathbb C} V_1$
and $\dim_{\mathbb C} V_1= 2 + a_1+2b_1= \rho_1+\rho_2$
since $V_1$ is a Jordan triple with characteristic
$(2, a_1, b_1)$, and $\rho_{\fg}= 1+\rho_1+\rho_2$.
The fact about the Cayley transform 
is well-known; see e.g. \cite{Sch}, \cite[Proposition 10.6(3)]{Loos-bsd}. 
\end{proof}

We note that the evaluation 
of character $\chi$  on $H_0, H_1, H_2$ is given by 
\begin{equation}
\label{chi-H-E}  
\chi(H_0) =2\chi(Z) =2i, \, 
\chi(H_1) =\chi(H_2)=-i, 
\end{equation}
since $H_0e = 2Ze=2ie$, 
$H_je = i(D(v_j,  v_j) -D(e,  e))e= -ie$. 
We observe and keep in mind that the geometrically natural choice of 
the complex structure of $S_1=K/L_0=K_1/L_1$
results in some discrepancy: For 
$H_o=iD(e,  e), \ad (-iH_0)=\ad (D(e,  e))$, 
has 
{\it non-negative} eigenvalues $2, 1, 0$ on $\fp^+=V 
=V_2+V_1+V_0$, 
whereas it has {\it negative} eigenvalue $-1$ on $\fq^+$, 
\begin{equation}
  \label{Dee-Dve}
[D(e,  e), D(v,  e)] = - D(v,  e), \quad D(v,  e)\in \fq^+.  
\end{equation}

\subsection{Cartan-Helgason Theorem for $K/L_0$}

Let $L^2(K, L_0, \chi_l)$
be the  $L^2$-space 
of sections of the homogeneous
line bundle $K\times_{ (L_0, \chi_l^{-1}   )  }\times \mathbb C $
defined by the 
$\chi_l^{-1}$ of $L_0$. The space $L^2(K, L_0, \chi_l)$ consists of $f\in L^2(K)$ such that
\begin{equation}
  \label{def-L2-chi}
f(gh) =\chi^l(h)f(g), \quad g\in K, h\in L_0, \, he=\chi(h)e. 
\end{equation}
It follows immediately from the definitions of $L_0$ and $\chi$
that
\begin{equation}
  \label{chi-def}
L^2(K/L)=\sum_{l=-\infty}^{\infty}L^2(K, L_0, \chi_l). 
\end{equation}
 under the regular left action  of $K$.

We shall treat extensively functions on $K$
that are transforming under $L_0$
as in (\ref{def-L2-chi} ) and it is convenient
to give the following

\begin{defi+}
An element $f\in L^2(K)$ is called $(l_1, l_2)$-spherical if
  $$
  f(h_1kh_2) =
  \chi_{l_1} (h_1) \chi_{l_2}(h_2) 
f(k), \quad h_1, h_2\in L_0.
  $$
\end{defi+}

\begin{lemm+} \label{C-H}
  \begin{enumerate}
  \item
    Let $\fg\ne \fsu(r+b, r)$, $r>1$.
  The space    $L^2(K, L_0, \chi_l)$
  is decomposed as 
  a sum of irreducible representations of $K$, 
  $$
  L^2(K, L_0,  \chi_l) 
  =\sum_{\mu}W_{\mu, l}, 
  $$
  where each $W_{\mu, l}$
  has highest weight  given by 
\begin{equation}
  \label{mu-m}
  {l\alpha_0+\mu,  \,\, \mu= (\mu_1, \mu_2)=\mu_1\alpha_1 + \mu_2\alpha_2, 
    \, 
    \quad 
\mu_1\ge \mu_2\ge |l|, 
\quad \mu_1= \mu_2 = l, \mod 2 }. 
\end{equation}
  Moreover each space 
  $W_{\mu, l}$
  contains a unique vector $(l, l)$-spherical element
  $\phi_{\mu, l} $.
\item     Let $\fg= \fsu(r+b, r)$, $r>1$.
  The space    $L^2(K, L_0, \chi_l)$
  is decomposed  as above with
$$\mu_1,  \mu_2\ge |l|
\quad \mu_1= \mu_2 = l, \mod 2.
$$
The highest weight vector in $W_{\mu, l}$
can be chosen as
$$
f(x) = z_1^{p} {\overline{z_{r+b}}}^{q} w_1^{p'}{\overline{w_{r}}}^{q'}, 
\quad x=zw^\ast\in S \subset  M_{r+b, r}(\mathbb C), 
$$
where we have written a rank-one projection $x\in S\subset  M_{r+b, r}(\mathbb C)$
as $x=zw^\ast, z\in \mathbb C^{r+b}, w\in \mathbb C^{r}$, $\Vert z\Vert =
\Vert w\Vert =1$, and where $(p, q, p', q')$ are subject to the
condition
$$
\mu_1=p+q, \mu_2=p'+q',  l=p-q=p'-q'.
$$
The spherical polynomial   $\phi_{\mu, l}$ in this case is
$\phi_{\mu_1, -l}(zw^\ast) =\phi_{\mu_1, -l}(z)  \phi_{\mu_2,l}(z)$
where $\phi_{a, -l}$
is $(l, l)$-spherical polynomial
on $\mathbb P(\mathbb C^{n})$.
\end{enumerate}
\end{lemm+}

\begin{proof}
  \begin{enumerate}
  \item   The statement for 
  the decomposition of $L^2(K, L_0,  \chi_l)$
as   representation   of the semi-simple group $K_1$ is 
in \cite[Theorem 7.2]{Sch}; our $(\fk_1^\ast, \fl_1)$ corresponds
to $(\fg, \fk)$ there.
More precisely our character $\chi_l$ 
is precisely the same as 
$\chi_{-l}$ in 
\cite{Sch}.
The character 
$\chi_{-l}$  on $K$ in \cite{Sch} 
for the Hermitian symmetric space $D=G/K\subset \mathbb C^d=\fp^+$ 
is defined  
by  
$$ 
X\in \fk^{\mathbb C} \mapsto \frac l{\tr \Ad D(e, e)|_{\fp^+}}\tr \Ad X|_{\fp^+};
$$ 
equivalently 
it is determined  \cite[(5.1)]{Sch} by $\chi_{-l}:  D(e, e)\mapsto l$.
For our symmetric pair $(\fk_1^\ast, \fl_1)$
the corresponding $D(e, e)$
is  $H_1 =D(v_1,  {v_1}) -D(e, e)$
described in Lemma  \ref{k-roots} and
$H_1 e= -l$ and thus $\chi_l(H_1)=-l$. 
(Alternatively we can also prove this  by computing
using the duality relation in Appendix \ref{app-b}.)
The
results in \cite{Sch} then 
determine the highest weights of $W_{\mu, l}$
on $\fh_{i\fq}=\mathbb R(iE_1) +\mathbb
R(iE_2)$ as $\mu_1\alpha_1 +\mu_2\alpha_2$
in our statement.

Finally it  is trivial to find the weight of $W_{\mu, l}$ on the central element $Z$.
The right action of  $\exp (sZ)$ on 
  $\phi\in W_{\mu, l}$
is, using $Ze=ie$, 
\begin{equation}
  \label{Z-action}
  \begin{split}
\pi_{\nu}(\exp (sZ))\phi(h) 
 & =
  \phi( \exp (-sZ) h) =
  \phi(h \exp (-sZ) ) 
  =    \exp(-isl) \phi(h).  
 \end{split}
 \end{equation}
  Thus 
$\pi_\nu(Z)\phi 
  = - il\phi, 
  $  and
  $\pi_\nu(iZ)\phi 
  = l\phi $.
  Thus the highest weight of
   $W_{\mu, l}$ is
  $l\alpha_0+\mu_1\alpha_1 +\mu_2\alpha_2$.

\item This part is well-known; see e.g.  \cite{HT, JW}.
\end{enumerate}
\end{proof}

Altogether we have now $I(\nu)$ is
  $$
  I(\nu)=L^2(K/L)
  =\sum_{l=-\infty}^\infty\sum_{\mu}W_{\mu, l}
  $$
with $\mu$ being specified above.

  \begin{rema+} The exact formulas for the highest weight vectors 
  above in  the case  $\fg= \fsu(r+b, r)$ are not needed in our paper 
  However it is possible to prove  Theorem \ref{main}
  below by using the weight vectors instead of spherical vectors; see \cite{HT}.
Note also that the parametrization of the spherical
  polynomials  on $\mathbb P(\mathbb C^{n})$ generated by 
the $(p, q)$-spherical harmonic polynomial $  z_1^{p} {\overline{z_{n}}}^{q}
$ as  $\phi_{\mu_1, -l}$, $\mu_1=p+q, l=p-q$,
is due to the our geometric definition  of character $\chi$. Recall that due to (\ref{chi-H-E})
$\phi_{\mu, l}$ satisfies
\begin{equation}
  \label{chi-H-j}
\phi_{\mu, l} (k e^{it H_j}) = e^{-ilt}  \phi_{\mu, l} (k), \quad, j=1, 2, 
\end{equation}
The same parametrization is used in  Appendix \ref{app-a}.
\end{rema+}

\subsection{Harish-Chandra $C$-function and
  expansion of the 
  $(l, l)$-spherical
polynomials
}

A major technical step in the proof
of Theorem    \ref{main}  below is to use the Harish-Chandra $c$-function 
to compute certain expansions and differentiations involving 
the spherical polynomials $\phi_{\mu, l}$. 
We recall that the  spherical polynomial 
$\phi_{\mu, l}(h) =\phi_{\mu, l}(h) $ on  $K/L=K_1/L_1$
 is a special case  
of the  Harish-Chandra spherical  
function $\Phi_{\lambda, l}(h)$
for the non-compact pair $(\fk_1^\ast, \fl_1)$
corresponding to the one dimensional character $\chi_l$
of $L_1\subset L_0$;  see \cite{Sch, Shi}. 

 More  
precisely \cite{He3, Sch, Shi}
\begin{equation}
  \label{phi-Phi}
  \phi_{\mu, l}(h)
  =\Phi_{-i(\mu + \rho  ),
    l}(h),  h\in K_1.
\end{equation}
where $\rho=\rho_{\fk_1^\ast}$.
The spherical function $\Phi_{\lambda, l}$ is
invariant with respect to the Weyl group
$W(\fk^\ast_1, \fh_{i\fq})$ of the root system
$R(\fk^\ast_1, \fh_{i\fq})$
in (\ref{root-sys-q}).
Eventually we shall replace $\lambda$ by $-i(\mu + \rho)$
and use {\it Weyl  group symmetry in $\mu + \rho$}.
 The leading term of 
$\Phi_{\lambda, l}(h)$ 
is given by the limit formula 
\begin{equation}
  \label{Phi-HC-c}
\lim_{H\to \infty} 
e^{-(i\lambda -\rho)(H)}
\Phi_{\lambda, l}(\exp(H))= C(\lambda, l), 
\end{equation}
for $H$ in the positive Weyl Chamber of 
the root system (\ref{root-sys-q}), i.e, for $H=x_1 (iE_1) +
x_2(iE_2)\in \fh_{iq}, x_1 >x_2 >0$, and
for $\re(i\lambda)= y_1\alpha_1 + y_2\alpha_2, y_1>y_2>0$; see 
\cite[Ch.~IV, Theorem 6.14; Ch.~V, Section 4]{He3}, \cite[Theorem 3.6]{Shi}.
In particular
\begin{equation}
  \label{Phi-HC-c-var}
  \phi_{\mu, l}(\exp(H))=
  \Phi_{\lambda, l}(\exp(H))= C(\lambda, l) e^{(i\lambda
    -\rho_{})(H)} + L.O.T., \, H\in
  \fh_{\fq}^{\mathbb C}
\end{equation}
as an expansion of trigonometric polynomial on the complexification 
$\exp(\fh_{\fq}^{\mathbb C})\cdot [e]$
of the
real torus $\exp(\fh_{\fq}) \cdot [e]\subset K/L_0
=\mathbb P(K/L), [e]=\mathbb C e$, with lower order terms (L.O.T.) being trigonometric polynomials
of lower order in the sense defined by the Weyl Chalmers.
Here 
 $C(\lambda, l)=C(\lambda, -l)$
 is the Harish-Chandra $C$-function, 
 and in our case 
it is given by 
$$ 
C(\lambda, l) 
=c_0 
\prod_{\epsilon=\pm 1}\frac 
{\Gamma(\frac i2 (\lambda_1 +\epsilon\lambda_2)) 
}
{
\Gamma(\frac 12 a_1 + i(\lambda_1 +\epsilon \lambda_2))  
}
\prod_{j=1,2}
\frac 
{2^{-i\lambda_j} \Gamma (i\lambda_j)}
{
  \Gamma(
  \frac 12 
  (b_1 +1 + i\lambda_j +l) 
  ) 
  \Gamma(
  \frac 12 
  (b_1 +1 + i\lambda_j -l) 
)  
}
$$ 
for $\lambda=\lambda_1 \alpha_1 + \lambda_2\alpha_2$, 
where $c_0$ is normalized so that $C(-i\rho, 0)=1$
for the Harish-Chandra $C$-function $C(-i\rho, 0)$
with  trivial line bundle, $l=0$; 
see loc.~cit.. We observe also that the $C$-function 
is positive $C(\lambda, l) $ for $\lambda =-i(\mu +\rho)$.

In particular we shall need
the spherical polynomials $\phi_{(1, 1), \pm 1}(k)$ for $\mu=(1, 1)$.
The corresponding representations space of $K$
is $\fp^{+}=V$ or $\fp^{-}=\bar V$, and the 
spherical polynomial is  the matrix coefficient
\begin{equation}
  \label{mat-coe}
  \phi_{(1, 1),  1}(k) 
  =\langle k e, e\rangle, \quad 
\phi_{(1, 1),  -1}(k) =
\langle e, k e\rangle.
\end{equation}
Indeed  the space  $\fp^{+}=V$
is a representation 
of $\fk_1^{\mathbb C}$  of highest weight 
$\alpha_1 +\alpha_2$
and representation of $\fk^{\mathbb C}$  of highest weight
$\alpha_0+\alpha_1 +\alpha_2$ since $Z$ acts as $i$,
the corresponding highest weight vector  is
\begin{equation}
  \label{hw-H}
v_{0}= \frac 12 \left((v_1 -i e) + (v_2 + iw)\right). 
\end{equation}
In other words, recalling 
that $e$ is the root vector of 
the Harish-Chandra strongly orthogonal root 
$\gamma_1$, we see that 
$\alpha_1 +\alpha_2$ is conjugated to $\gamma_1$.





\section{Lie algebra $\fg$-action on 
  $I(\nu)=Ind_P^G( \nu)  $
}
\label{Liealg-act}
We compute the Lie algebra
action of $\fg^{\mathbb C}$ on
$I(\nu)  $.
For that purpose 
we denote the right differentiation
of  Lie algebra elements $X\in \fg^{\mathbb C}$
on  functions $f$ on $G$ by $Xf$,
\begin{equation}
  \label{right-diff}
Xf(g)=\frac{d}{ds}f(g\exp(sX))|_{s=0}.  
\end{equation}
Then $X$ commutes with the left regular action
\begin{equation}
  \label{chain-rule-0}
 X(f)(hx)= X\(f(h\cdot )\) (x),
\end{equation}  
and intertwines the right  action  $f(x)\to f(xh)=f_h(x)$ as
\begin{equation}
     \label{chain-rule}
     (Xf)_h(x)=   (Xf)(xh)
     = \( (\text{ad}\, h (X) )  f_h\)(x).
   \end{equation}

First it follows
from    (\ref{xi-H0}) and  (\ref{def-chi}) that $H_0 e=2ie$, $\chi_l(\exp(tH_0))=e^{2ilt}$, 
$\chi_l(H_0)=2il$.
Thus any element
 $f\in L^2(K, L_0,  \chi_l)$ is an eigenfunction
of the differentiation by $H_0$,
\begin{equation}
  \label{f-under-H_0-var}
H_0f = 2il f.
\end{equation}

\begin{theo+}
  \label{main}
  Let   $\fg$ be a simple Hermitian Lie algebra of rank $r\ge 2$
  and $\fg\ne \fsp(r, \mathbb R)$.
The action of $  \pi_\nu(\xi)$ on 
$\phi_{\mu, l}$ is given by
\begin{equation*}
\begin{split}  
&\quad 2^3\,  \pi_\nu(\xi) 
\phi_{\mu, l}
\\
&
=
\sum_{\sig=(\sig_1, \sig_2)=(\pm 1, \pm 1)}
\(\nu  +{\sig_1}(\mu_1+\rho_1) 
+\sigma_2 (\mu_2+\rho_2) 
 -(\rho_1+\rho_2)     \) 
 \\
 &\times \(
 c_{\mu, l} (\mu+\sig, l+1)
\phi_{\mu +\sig, l+1}
+c_{\mu, l}(\mu+\sig, l-1)
\phi_{\mu +\sig, l-1}\), 
\end{split}
\end{equation*}
where  the coefficients 
$
c_{\mu, l} (\mu +\sig, l\pm1)
$
are  given by 
\begin{equation}
  \label{hc-c-frac}
c_{\mu, l} (\mu +\sig, l+1) 
=
\frac 
{C 
  \left(
    -i(s_{\sig}(\mu +\rho)) 
,l \right) 
}
{ C 
  \left(
    -i 
    (\sigma_1\alpha_1+\sigma_2\alpha_2+
    s_{\sig}(\mu +\rho)),   l+1 
    \right) 
  }, 
\end{equation}
\begin{equation}
  \label{hc-c-frac-2}
c_{\mu, l} (\mu+\sig, l-1) 
  =
 \frac 
{C \left(
    -i(s_{\sig}(\mu +\rho)) 
,l  \right) 
}
{ C 
  \left(
    -i 
    (\sigma_1\alpha_1+\sigma_2\alpha_2+
    s_{\sig}(\mu +\rho)),   l-1 
    \right) 
}.
\end{equation}
and
$s_{\sig}$ is in the Weyl group
$W$ of the root system
  (\ref{root-sys-q})
such that $s_{\sigma}(\alpha_1+\alpha_2)=
\sigma_1\alpha_1+\sigma_2\alpha_2$.  Moreover
all the coefficients are positive.
\end{theo+}

It is understood here that  the  term
$ \phi_{\mu +\sigma, l\pm 1}= \phi_{(\mu_1 +\sigma_1,  \mu_2+\sigma_2), l\pm 1}$
will not appear in the RHS
if $(\mu_1 +\sigma_1)\alpha_1   +(\mu_2+\sigma_2)\alpha_2$
is not one of the highest weights specified in Lemma \ref{C-H}.

\begin{rema+}
  \label{rema+-}
  It is remarkable that all the coefficients 
of $\phi_{\mu +\sig, l\pm 1}$
have a rather uniform 
formula.  Actually it is relatively easy to find
the coefficient of the leading term
$\phi_{\mu + (1, 1), l\pm 1}$ and the other coefficients
can be obtained from the Weyl group symmetry
and by unitarity of $\pi_{\nu}$
for $\nu=\rho_{\fg} +ix, x\in \mathbb R$.
We shall find all the coefficients independent
of the unitarity, and we prove some recursion formulas
for spherical polynomials which might be of independent interests \cite{Z-tams}.

\end{rema+}

\begin{proof}
 We claim first that for any $X\in \fp^{\mathbb C}$, 
$$
\pi_\nu(X) 
W_{\mu, l}
\subseteq
\sum_{\sig_1, \sig_2=\pm 1)} W_{\mu +\sig, l\pm 1}.
$$
This follows by considering the tensor product 
$\fp^{\mathbb C}\otimes W_{\mu, l}$
as representation of $K$. 
Indeed let $\fg\ne \fsu(r+b, r)$.
The
adjoint action of
the central element $Z\in \fk$ on
$\fp^{\pm}$ is $\pm i$, and its right action on 
$W_{\mu, l}$ is $il$.
Let $X\in \fp^+$, then
$\pi_\nu(X) 
W_{\mu, l}$
is of weight
$i(l+1)$ under $Z$ for any $X\in \fp^+$.
It is also a classical fact 
that the highest weights in 
the tensor product  decomposition of
$W_{\mu, l}\otimes \fp^{\pm}$
under $\fk^{\mathbb C}$
are of the form $\mu+l\alpha_0+\nu'$ where $\nu'$ is a
weight 
appearing in $\fp^{+}$.
The space $\fp^{+}$
is of highest weight $
(1, 1)=\alpha_1+
\alpha_2$ under the Cartan subalgebra
$(\fh_{i\fq})^{\mathbb C}=\mathbb CE_1 +\mathbb CiE_2$
of $\fh^{\mathbb C}$
and the only weights
in $\fp^{\pm}$ of the form $c_1\alpha_1 +
c_2\alpha_2$ 
are  $\sigma_1\alpha_1+
\sigma_2\alpha_2,
\sigma_1\alpha_1,  \sigma_2\alpha_2$. However 
  by the Cartan-Helgason theorem, Lemma \ref{C-H},  we see 
  that  $\sigma_1\alpha_1,  \sigma_2\alpha_2$
 are not eligible.
Thus $\pi_\nu(X) 
W_{\mu, l}$ is of the claimed form. 

When $\fg=\fsu(r+b, r)$
the Weyl group
for the root system of $(\fk^\ast, \fh_{i\fq})$
is $(\mathbb Z_2)^2$ consisting of only sign changes
instead of  all signed permutation $(\mathbb Z_2)^2\rtimes S_2$, but all the relevant
weights $\sigma_1\alpha_1+
\sigma_2\alpha_2$ are still in the orbit of the Weyl group
$(\mathbb Z_2)^2$ so  the arguments are valid for
$\fg=\fsu(r+b, r)$ as well.

Next, the element  $\xi =\xi_e$ is
  invariant under $L\subset K$, thus
$    \pi_\nu(\xi) 
\phi_{\mu, l}$
is  a sum of the $L$-invariant
vectors
in $\sum_{\sig_1, \sig_2=\pm 1} W_{\mu +\sig, l\pm 1}$,
and is further by Lemma  \ref{C-H}
a linear combination of
$\phi_{\mu +\sig, l\pm 1}$. The rest
of the proof is to determine the  coefficients. 

Notice that
each function 
in the linear combination
 is determined by its restriction on
  the complexification 
 $\exp(\fh_{\fq}^{\mathbb C})\cdot [e]$
 once the line parameter $l$ is given, so
 it is enough to find the expansion restricted on the torus
(after the differentiations) as the line bundle parameters of each term in the expansion
 are already fixed.

 We have 
\begin{equation}
  \begin{split}
\pi_\nu(\xi) 
\phi_{\mu, l} 
(k )
&=\frac{d}{ds}
\phi_{\mu, l}
(
\exp(-s\xi)
k
)
|_{s=0}=
\frac{d}{ds}
\phi_{\mu, l}
(kk^{-1}
\exp(-s\xi)
k
)
|_{s=0}
\\
&=\frac{d}{ds}
\phi_{\mu, l}
(
k
\exp(-s\ad(
k^{-1}
) \xi))|_{s=0}
=- \( (\ad (k^{-1}) \xi) \phi_{\mu, l} \)(k), \, k\in K,
\end{split}
\end{equation}
where $\( (\ad (k^{-1}) \xi) \phi_{\mu, l} \)(k)$  is
right differentiation of the Lie algebra valued 
{\it vector field} $-\ad (k^{-1}) \xi$
on  $\phi_{\mu, l}$  evaluated at $k\in K$.
The element
$\phi_{\mu, l}$
is in the induced representation any differentiation
of $\phi_{\mu, l}$  along the Lie algebra $\fm + \fn$ is zero, and 
we need formulas for  $\ad (k^{-1}) \xi=\ad (k^{-1}) \xi_e=
\xi_{k^{-1} e}
$ mod
$\fm + \fn$.

\begin{lemm+}
  \label{mod}
Let 
    $V= V_2+ V_1+ V_0=\mathbb Ce 
    + V_1+ V_0$ be the 
    Peirce decomposition with respect to  the minimal 
tripotent $e$ and 
    $P_2, P_1, P_0$  
the corresponding projections. 
Any element     $\xi_{u}\in \fp$
has the following decomposition according to 
(\ref{nman}), $\mod \fm + \fn$, 
\begin{equation}
\label{k-xi}  
\xi_{u} = \re 
\langle  u , e\rangle \xi 
+ \im 
\langle  u, e\rangle H_0  + D(P_1(u),  e) -
D(e, {P_1(u)}). 
\end{equation}
\end{lemm+}  
\begin{proof}
Write $u = P_2u + P_1u + P_0u 
  =\langle u, e\rangle e + u_1 + u_0
  =\re
  \langle u, e 
  \rangle e
  +  i \im \langle u, e  
  \rangle   e + u_1 + u_0
  $. 
Then
  $$
\xi_{\langle u, e\rangle e}
=\re \langle u, e\rangle \xi_{e}
+ \im \langle u, e\rangle \xi_{ie}, 
$$
with 
$$
\xi_{ie}=(\xi_{ie} - H_0)+ H_0 = H_0, \, \mod \fn 
$$
by (\ref{n2}). In view of (\ref{n1}) we have 
$$
\xi_{u_1}
=\(\xi_{u_1}  + D(e,  u_1) -D(u_1,  e)\) +
\(D(u_1,  e) -D(e,  u_1)\) =
D(u_1,  e) -D(e,  u_1), \mod \fn, 
$$
and  $\xi_{u_0}\in \fm$. This proves  (\ref{k-xi}). 
\end{proof}

Using Lemma \ref{mod} and the
formula  (\ref{mat-coe})
we see that $
\ad (k^{-1}) \xi $, mod $\fm^{\mathbb C}+
\fn^{\mathbb C}$, is
\begin{equation}
  \label{pear-1}
  \begin{split} 
  &  \quad \ad (k^{-1}) \xi
    \\
    &=\re\langle k^{-1}e, e\rangle \xi
+\im\langle k^{-1}e, e\rangle H_0
+ D(P_1(k^{-1}e), e) -D(e, P_1(k^{-1}e))\\
&=\re \langle k e, e\rangle
\xi
-\im\langle k e, e
\rangle H_0
+ D(P_1(k^{-1}e), e) -D(e, P_1(k^{-1}e)).
\end{split}
\end{equation}
Hence
\begin{equation}
  \label{apple-1}
 \pi_\nu(\xi) \phi_{\mu, l} (k) = I + II + III
\end{equation}
with
$$
I=-\re \langle k^{-1}e, e\rangle \(\xi  \phi_{\mu, l}\)(k),$$
$$
II=\im \langle k e, e\rangle    \({H_0}  \phi_{\mu, l}\)(k)
$$
and
$$III=\(
\[D(e, P_1(k^{-1} e)) 
-D(P_1(k^{-1} e), e)
\]\phi_{\mu, l}\)(k).
$$
Using the definition  (\ref{def-ind})
of the induced representation 
we have the right differentiation of $\xi$
on any element $f\in L^2(K/L)=Ind_P^G(\nu)$
has eigenvalue $-\nu$,
$\xi  \phi_{\mu, l} =-\nu\phi_{\mu, l}$, and
 the first term is
 \begin{equation*}
   \begin{split}
 I    &=    -\re \langle k^{-1}e, e\rangle
     \(\xi  \phi_{\mu, l}\)(k)
\\
&=     \nu    \re \langle k^{-1}e, e\rangle
\phi_{\mu, l}(k)
\\
&=\frac{\nu}{2} \langle k e, e\rangle 
\phi_{\mu, l}(k) 
+\frac{\nu}{2} \langle e, k e\rangle
\phi_{\mu, l}(k)
\\
&=I^+ + I^-
   \end{split}
 \end{equation*}
 with $$
I^+= \frac{\nu}{2} \langle k e, e\rangle 
\phi_{\mu, l}(k) =\frac{\nu}{2} \phi_{(1, 1), 1}(k) 
\phi_{\mu, l}(k)
$$
and
$$I^-=\frac{\nu}{2} \phi_{(1, 1), -1}(k) 
\phi_{\mu, l}(k)  .
$$

The second term $II$, in view of   (\ref{f-under-H_0-var}), is 
\begin{equation*}
   \begin{split}
II&
=\im \langle k e, e\rangle (2il \phi_{\mu, l})(k) 
\\
&=l\langle k e, e\rangle \phi_{\mu, l}(k) 
-l\langle e, k e\rangle\phi_{\mu, l}(k)\\
&=l\phi_{(1, 1), 1}(k) 
\phi_{\mu, l}(k) 
-l\phi_{(1, 1), -1}(k)\phi_{\mu, l}(k) 
=II^+
+ II^-.
\end{split}. 
\end{equation*}

We proceed to find recursion formulas
for $\phi_{(1, 1), 1}\phi_{\mu, l}$.
For that purpose we find explicit coordinates
for the  complex torus
 $\exp(\fh_{\fq}^{\mathbb C})\cdot [e]$,
$\exp(\fh_{\fq}) e\subset S=K/L\subset V$.
We shall treat all functions as trigonometric functions
on the compact homogeneous space $K/L$,  the evaluation and
extension to the complexification
 $\exp(\fh_{\fq}^{\mathbb C})\cdot [e]$ being done differentiations.

\begin{lemm+}
If $k=\exp(x_1 E_1 + x_2 E_2), x_1, x_2\in \mathbb C$,  then 
$$
k^{-1} e =
\cos x_1   \cos x_2 e 
-(\sin x_1\cos x_2 v_1+\sin x_2\cos x_1 v_2)  + \sin x_1\sin x_2 w. 
$$
\end{lemm+}
\begin{proof}
  Using    (\ref{quadruple}) we find
  $$
E_1 e=  (D(v_1, e)- D(e, v_1))e = v_1, 
E_1^2 e=  E_1 v_1 =(D(v_1, e)-D(e, v_1)) v_1 =-e 
  $$
  and  generally 
  $$
    E_1^{2m}e = (-1)^me, \, 
  E_1^{2m+1}e = (-1)^m v_1.
  $$
  Therefore 
  $$
  e^{-x_1 E_1 }  e=
  \cos x_1 e  - \sin x_1 v_1, 
  $$
and also $ e^{-x_2 E_2  }  e=  \cos x_2 e -\sin x_2 v_2$. We compute further 
  $$
  E_2 v_1 = (D(v_2, e)-D(e, v_2))v_1 = w,\, 
  E_2^2 v_1= (D(v_2, e)-D(e, v_2))w = -v_1, 
  $$
  and in general 
  $$ 
    E_2^{2m} v_1 = (-1)^{m} v_1, 
  E_2^{2m+1}v_1  = (-1)^{m} w. 
  $$
  This implies that 
  $$
  e^{-x_2 E_2  } \cdot v_1= \cos x_2 v_1 -\sin x_2 w. 
  $$
  We have then 
  \begin{equation}
      \label{ke}
    \begin{split}
k^{-1} e&=e^{-x_1 E_1-x_2 E_2  } e =
e^{-x_2 E_2  } \cdot  (  \cos x_1 e -\sin x_1 v_1)\\
&  = 
  (\cos x_1 (  \cos x_2 e -\sin x_2 v_2) 
  +\sin x_1 (\cos x_2 v_1 -\sin x_2 w) 
  \\
  &=\cos x_1 \cos x_2 
e -(\sin x_1\cos x_2 v_1+
\sin x_2\cos x_1 v_2)  + \sin x_1\sin x_2 w.
  \end{split} 
 \end{equation}
\end{proof}

 \begin{lemm+}
   \label{lemma-2}
  The 
  following recursion formulas
  hold,
 \begin{equation}
\label{exp-1}
\phi_{(1, 1), 1}\phi_{\mu, l}
=\frac 14
\sum_{\sig_1, \sig_2=\pm 1}
c_{\mu,    l}(\mu+\sig, l+1)
\phi_{\mu +\sig, l+1},
 \end{equation}
 \begin{equation}
\label{exp-2}
\phi_{(1, 1), -1}\phi_{\mu, l}
=\frac 14
\sum_{\sig_1, \sig_2=\pm 1}
c_{\mu,  l}(\mu+\sig, l-1)
\phi_{\mu +\sig, l-1},
 \end{equation}
 where $ c_{\mu,  l}(\mu+\sig, l\pm 1)$ are given
 in Theorem 3.1.
\end{lemm+}

\begin{proof} Observe again  that
  by general tensor product
  arguments the product $\phi_{(1, 1), \pm 1}\phi_{\mu, l}$
  is a sum of
  $\phi_{\mu +\sigma, l-1}$, $\sigma=(\pm 1, \pm 1)$.
  We use the idea in \cite{Vretare}
  by considering the leading term
of $\phi_{(1, 1), 1}
\Phi_{\lambda, l}$ by using
 the Harish-Chandra $c$-function;   see also \cite{Z-matann}.
   We recall
     (\ref{phi-Phi}) and
   consider the
   expansion
   \begin{equation}
  \label{a-phi}
 \phi_{(1, 1), 1}
  \Phi_{\lambda, l}
= \sum_{\sig_1, \sig_2=\pm 1 } A_{\lambda -i\sig,     l+1}
\Phi_{\lambda +\sig, l+1}  + L.O.T..
\end{equation}
(Presumably L.O.T. will  not appear for {\it general} $\lambda$
but it will not concern us here.)
We let $h=\exp(H), H=x_1 (iE_1) + x_2(iE_2)$.
First it is clear from (\ref{mat-coe})
and (\ref{ke})
that
   \begin{equation}
  \label{elem-phi}
\phi_{(1, 1), 1} (h)
=\langle he, e\rangle
=\cosh x_1 \cosh x_2
= \frac 14 (e^{x_1 } + e^{-x_1 }) (e^{x_2 } + e^{-x_2 })
= \frac 14 e^{x_1 +x_2} + L.O.T..
   \end{equation}
The coefficient $\frac 14$ can also
be obtained using
the general formula
(\ref{Phi-HC-c});
indeed the evaluation of Harish-Chandra $C$-function  is
$$
C((1, 1), 1) = C(\alpha_1+\alpha_2, 1) =
\frac 14.
$$
Using the limit formulas (\ref{Phi-HC-c}) again we see that
the coefficient
$A_{\lambda -i (1, 1), l+1}$
of the leading term
$\Phi_{\lambda -i(1, 1), l+1}$
is
$$
A_{\lambda -i(1, 1), l+1}
=\frac 14
\frac{C(\lambda, l)}{
  C(\lambda -i(\alpha_1+\alpha_2),
  l+1) 
  }.
  $$
  We  use then the Weyl group symmetry
  $\Phi_{\lambda, l}=
  \Phi_{s(\lambda), l}$, $s\in \mathbb Z_2\times
  \mathbb Z_2=\{\pm 1\}\times \{\pm 1\}$
  to get
$$
A_{\lambda -i\sig,     l+1}
=
\frac 14
\frac{C(s_{\sig}\lambda, l)}{
  C(s_{\sig}\lambda -i(\alpha_1+\alpha_2),
  l+1) 
  }.
$$
The coefficients in the expansion (\ref{exp-1})
are then $
=A_{\lambda -i\sig,  l+1}$ with $\lambda =-i(\mu +\rho)$. The second
expansion (\ref{exp-2}) is proved by the same method.

\end{proof}

We can now apply the lemma to both terms $I$ and $II$,
\begin{equation}
  \label{I-upd-1}
I^+= \frac{\nu}{2}
\phi_{(1, 1), 1}(k) 
\phi_{\mu, l}(k)=\frac{\nu}{2^3}
\sum_{\sig_1, \sig_2=\pm 1}c_{\mu, l}(\mu+\sigma,    l+ 1) 
\phi_{\mu +\sig, l+1}(k), 
\end{equation}
\begin{equation}
  \label{I-upd-2}
I^{-}=\frac{\nu}{2} \phi_{(1, 1), -1}
\phi_{\mu, l}(k) 
=\frac{\nu}{2^3}
\sum_{\sigma_1, \sigma_2=\pm 1}c_{\mu, l} (\mu+\sig,    l- 1) 
\phi_{\mu +\sig, l-1}(k),
\end{equation}
\begin{equation}
  \label{II-upd-1}
 II^+=l\phi_{(1, 1), 1}(k)
 \phi_{\mu, l}(k)
= 
 \frac l 4\sum_{\sig_1, \sig_2=\pm 1} 
  c_{\mu, l}(\mu+\sig,    l+ 1)
\phi_{\mu +\sig, l+1},
 \end{equation}
\begin{equation}
  \label{II-upd-2}
 II^-=   -l\phi_{(1, 1), -1}(k) \phi_{\mu, l}(k)
= - \frac l 4\sum_{\sig_1, \sig_2=\pm 1} 
  c_{\mu, l}(\mu +\sig,    l- 1)
\phi_{\mu +\sig, l-1}.
\end{equation}

The third term $III$ 
is $$
III=\(D(e, P_1(k^{-1} e))\phi_{\mu, l}\)(k) 
-\(D(P_1(k^{-1} e), e)\phi_{\mu, l}\)(k) 
: =III^+
+ III^-
$$

\begin{lemm+}\label{diff} We have
  the following recurrence formula
  for the right differentiations
  of the vector fields $D(e, P_1(k^{-1}e)) $
  and 
  $-D(P_1(k^{-1}e), e) $
  on  $\phi_{\mu, l} $,
  \begin{equation}
    \label{diff-l+1}
    III^+
    =  \(
  D(
  e, P_1(k^{-1}e) 
  ) 
 \phi_{\mu, l}\) 
 (k)
 =
 \sum_{\sigma_1, \sigma_2=\pm 1} b_{\mu +\sig, l+1}
  \phi_{\mu +\sig, l+1}(k),    
  \end{equation}
  \begin{equation}
    \label{diff-l-1}
III^-=  -\(
  D(P_1(k^{-1}e), e)
    \phi_{\mu, l}  \)
  (k)
  =\sum_{\sigma_1, \sigma_2=\pm 1}  b_{\mu +\sig, l-1} 
  \phi_{\mu +\sig, l-1}(k),
\end{equation}
where the coefficients $b_{\mu +\sig, l\pm 1}$
are given by
$$
b_{\mu +\sig, l+1}
=
\frac 1{2^3}
\(
\sigma_1 (\mu_1+\rho_1) +\sigma_2 (\mu_2+\rho_2) -
(\rho_1+\rho_2) -2l)
\) 
c_{\mu, l}(\mu+\sig,    l+1),
$$
$$
b_{\mu +\sigma, l-1}
=
\frac 1{2^3}
\(
\sigma_1 (\mu_1+\rho_1) +\sigma_2 (\mu_2+\rho_2)
-(\rho_1+\rho_2) +2l
\) 
c_{\mu, l}(\mu+\sig,    l-1),
$$
\end{lemm+}


\begin{proof}
  Denote $X$ the vector field
  \begin{equation}
    \label{vect-X}
   X(k)=
D(e, P_1(k^{-1}e))  
  \end{equation}
  acting on functions on $K$ by right differentiation,  $f\to
  (X(k)f)(k)$. With slightly abuse
  of notation we abbreviate it sometimes as $(Xf)(k)$, 
  $III^+=(X(k) \phi_{\mu, l})(k) =(X\phi_{\mu, l})(k)$.
  We prove first that  $X \phi_{\mu, l}$
  is $(l+1, l+1)$-spherical.  Some care has to be taken as $X$ is 
  vector field taking values in the Lie algebra of $\fk^{\mathbb C}$;
  the transformation rule of  $X \phi_{\mu, l}$ under the center of $K$ in $L_0$
  is easily checked but  we have to prove
  it for all $L_0$.
The space  $ V_1$ is  invariant under the subgroup
  $L_0\subset K$, and 
  $$
  P_1((hk)^{-1}e)=  P_1(k^{-1}h^{-1}e)=
\chi(h)^{-1}  P_1(k^{-1} e),
\, h P_1((kh)^{-1}e)= P_1(k^{-1}e), h\in L_0.
$$
Also elements $h\in K$ act on $D(x, y)$ as Jordan triple
automorphisms  $\ad(h) D(u, v)= D(hu, hv)$, $D(u, v)$
is conjugate linear in $v$, and $\overline{\chi(h)} =\chi^{-1}(h), h\in L_0$;
elements  $k\in L$ act as Jordan triple isomorphism as $Le=e$
and $L_0$ as isomorphism up to the  character  $\chi$.
Thus the vector field $X(k)=D(e, P_1(k^{-1}e))  $
satisfies
  $$X(hk)  = {\chi(h)} X(k), \quad
  \ad(h) \(X  (kh)\)
  =\chi(h)  X(k), \quad h\in L_0.$$
  It follows 
  by  the chain rules  (\ref{chain-rule-0})
  and  (\ref{chain-rule}),
  \begin{equation*}
    \begin{split}
      ( X(hk) \phi_{\mu, l}\)(hk) &=
      \left(X (hk) \phi_{\mu, l}(h\cdot )\right)(k) \\
      &       = \left( \chi(h) X (k) \chi_l(h) \phi_{\mu, l}(\cdot)\right)(k)\\
      &= \chi_{l+1}(h) (X(k)\phi_{\mu, l})(k)
     \end{split}
  \end{equation*}
  and   \begin{equation*}
    \begin{split}      (X(kh)   \phi_{\mu, l})(kh)
&     =   \( \ad(h) X(kh)\)    \( \phi_{\mu, l}(\cdot h) \)(k)
    =\chi (h) \chi_{l}(h) 
    \( X(k)\phi_{\mu, l}\)(k)
    \\
    &=\chi_{l+1}(h) 
    \( X(k)\phi_{\mu, l}\)(k).
  \end{split}
\end{equation*}
    Thus $(X(k) \phi_{\mu, l}(k)$
  must be  of the form
 (\ref{diff-l+1}).
 To find the coefficients
 we consider the
 group $SL(2, \mathbb C)^2=SL(2, \mathbb C)\times SL(2, \mathbb C)$ with
 the Lie
 algebra
 $\fsl(2, \mathbb C)\oplus \fsl(2, \mathbb C)$
generated by   (\ref{H-12})
and the restriction
of $\phi_{\mu, l}$
on  $SL(2, \mathbb C)\times SL(2, \mathbb C)$
and its  expansion
in terms of spherical polynomial of
$SL(2, \mathbb C)$, namely we consider the branching
of the representation $(K^{\mathbb C}, W_{\mu, l})$
under $SL(2, \mathbb C)\times SL(2, \mathbb C)$.
The highest weight of
the representation $ W_{\mu, l}$
restricted
to
$SL(2)\times SL(2)$ is
$\mu=(\mu_1, \mu_2)=\mu_1\alpha_1+\mu_2\alpha_2$,
and thus the representation $\bigodot^{\mu_1}\mathbb C^2
\otimes \bigodot^{\mu_2}\mathbb C^2$  of $SL(2)\times SL(2)$
appears
in $(K^{\mathbb C}, \mu)$,  and all other representations
 are
 of the form $(\mu_1', \mu_2')$ with $\mu_1' <  \mu_1$.
 Let $\psi_{m, l}(g)$ be the $(l, l)$-spherical
 polynomials for the group $SL(2, \mathbb C)$
 in the Appendix \ref{app-a}, (\ref{psi-sl2}).
 Comparing the leading term of 
 $\phi_{\mu, l}(g_1, g_2) $
 and $\psi_{\mu_1, l}(g_1) 
\psi_{\mu_2, l}(g_2) $ we have then
 $$
\phi_{\mu, l}(g_1, g_2) 
= C(-i(\mu+\rho), l)
\psi_{\mu_1, l}(g_1) 
\psi_{\mu_2, l}(g_2) + L.O.T, \quad (g_1, g_2)\in SL(2, \mathbb C)^2.
$$
The vector field $X(k)=D(e, P_1(k^{-1}e))$
restricted to $k=(k_1, k_2)=
(\exp(x_1E_1),
\exp(x_2E_2))
\in SU(2)\times SU(2)$
is,  by Lemma \ref{mod},
of the form 
\begin{equation*}
  \begin{aligned}
X(k)
&=-\sin x_1 \cos x_2
D(e, v_1) 
-\sin x_2 \cos x_1 D(e, v_2)\\
&=-\sin x_1 \cos x_2  E_1^{-}
-\sin x_2 \cos x_1 E_2^{-}.    
  \end{aligned}
\end{equation*}
The vector field $X(k)$ takes values also in
the complexification of the Lie algebra
$\fsl(2, \mathbb C) +
\fsl(2, \mathbb C)$. Thus
the restriction of
$$
\(X\phi_{\mu, l}\)|_{SU(2)\times SU(2)}
=X \left(\phi_{\mu, l}|_{SU(2)\times SU(2)}\right).
$$
Moreover the Lie algebra differentiations
by $ D(e, v_1),  D(e, v_2)$ clearly 
preserve the degree. Thus we have
\begin{equation*}
  \begin{aligned}
\(X \phi_{\mu, l}\)(k_1, k_2) 
&= C(-i(\mu+\rho), l)
\(-\sin x_1 (E_1^-\psi_{\mu_1, l})(k_1)\) 
\(\cos x_2\psi_{\mu_2, l}(k_2)\) \\
&\quad 
+C(-i(\mu+\rho), l)
\(-\sin x_2 (E_2^- \psi_{\mu_2, l})(k_2)\) 
\(\cos x_1\psi_{\mu_1, l}(k_1)\)
+ L.O.T..
  \end{aligned}
\end{equation*}
We   use now Lemma \ref{A1},  (\ref{E-psi-expl}),
 and obtain
 \begin{equation}
   \label{III-key}
-\sin x_j (E_j^-\psi_{\mu_j, l})(k_j) 
=\frac 14 (\mu_j -l) \psi_{\mu_j+1, l+1}(k_j) + L.O.T.. , \quad j=1, 2
\quad
 \end{equation}
The leading term of $\cos x_2 \psi_{\mu_2, l}(k_2) $, $k_2=\exp(x_2E_2)$,
is clearly the same as
$\frac 12 \psi_{\mu_2+1, l}(k_2) $.
Thus
$$
\(-\sin x_1 (E_1^-\psi_{\mu_1, l})(k_1)\) 
\(\cos x_2\psi_{\mu_2, l}(k_2)\)
=\frac 18(\mu_1 -l)
\psi_{\mu_j+1, l+1}(k_j)\psi_{\mu_2+1, l}(k_2);
$$
similarly for $\(-\sin x_2 (E_2^- \psi_{\mu_2, l})(k_2)\) 
\(\cos x_1\psi_{\mu_1, l}(k_1)\)$. We have then
 \begin{equation*}
   \begin{split}
\(X \phi_{\mu, l}\)(k_1, k_2) 
&=
\frac 1{8} (\mu_1-l + \mu_2 -l) 
 C(-i(\mu+\rho), l) 
\psi_{\mu_1+1, l+1}(k_1) 
\psi_{\mu_2+1, 
  l+1}(k_2) 
+L.O.T. 
\\
&=\frac 1{8} (\mu_1+ \mu_2 -2l) 
 C(-i(\mu+\rho), l) 
\psi_{\mu_1+1, l+1}(k_1) 
\psi_{\mu_2+1, 
  l+1}(k_2) 
 +L.O.T..     
   \end{split}
 \end{equation*}

 On the other hand the leading term in  RHS
of    (\ref{diff-l+1}) is 
\begin{equation*}
  \begin{aligned}
&\quad b_{\mu + (1, 1), l+1}
\phi_{\mu + (1, 1), l+1}(k_1, k_2)\\
&=    b_{\mu +(1,1), l+1}
C(-i(\mu+ (1, 1)+\rho), l+1) 
\psi_{\mu_1+1, l+1}(k_1)\psi_{\mu_2+1, 
  l+1}(k_2) + L.O.T..
  \end{aligned}
\end{equation*}
It follows then that
$$
b_{\mu +(1, 1), l+1} =
\frac 1{8} (\mu_1 +\mu_2 -2l)c_{\mu, l}(\mu+(1, 1),    l+1)
$$
where $c_{\mu, l}(\mu+(1, 1),    l+1)$ is given in (\ref{hc-c-frac}).
This proves the formula for the leading
coefficient.

To find the other coefficients we write
$$
\mu_1 +\mu_2 -2l = (\mu_1+\rho_1) + (\mu_2+\rho_2) -(\rho_1 +\rho_2) -2l
$$
and  use the Weyl group symmetry
as in the proof of Lemma above
to get
$$
b_{\mu +\sig, l+1} =
\frac 1{8} \(\sigma_1 (\mu_1+\rho_1) +\sigma_2 (\mu_2+\rho_2) - (\rho_1+\rho_2)-2l)
\)
c_{\mu, l}(\mu+\sig,    l+1),
$$
for $\sigma_1, \sigma_2=\pm 1$.

To prove (\ref{diff-l-1}) we consider the vector field
$Y(k)=  -\(
  D(P_1(k^{-1}e), e)
    \phi_{\mu, l}  \)$ and its restriction to $SL(2, \mathbb C)^2$. We have
\begin{equation*}
  \begin{aligned}
Y=Y(k)
&= \sin x_1 \cos x_2  E_1^{+}
+\sin x_2 \cos x_1 E_2^{+}.    
  \end{aligned}
\end{equation*}
We
use
then (\ref{E+psi}) and \ref{E+psi-expl}
to find the leading term of the expansion
$Y\phi_{\mu, l}$ and obtain
all the coefficients
by  Weyl group symmetry.
  \end{proof}

  Hence
  \begin{equation*}
\begin{aligned}  
III^+ &=\(D(e, P_1(k^{-1} e))\phi_{\mu, l}\)(k) \\
&= \sum_{\sig_1, \sig_2=\pm 1}
\frac 1{2^3} \(
\sigma_1 (\mu_1+\rho_1) +\sigma_2 (\mu_2+\rho_2)
-(\rho_1+\rho_2) -2l
\) 
c_{\mu, l}(\mu+\sig,    l+1). 
  \phi_{\mu + \sig, l+1}(k),
 \end{aligned}  
\end{equation*}
\begin{equation*}
\begin{aligned}  
III^- & =-\(D(P_1(k^{-1} e), e)\phi_{\mu, l}\)(k) \\
&= \sum_{\sig_1, \sig_2=\pm 1}
\frac 1{2^3}
\(
\sigma_1 (\mu_1+\rho_1) +\sigma_2 (\mu_2+\rho_2) -
(\rho_1+\rho_2) +2l \) 
c_{\mu, l}(\mu+\sig,    l-1). 
\phi_{\mu +\sig, l-1}(k).
 \end{aligned}  
\end{equation*}

Altogether we find
$\pi_\nu(\xi)\phi_{\mu, l}
= (I^+ + II^+ + III^+) +
(I^- + II^- + III^-),
$ 
$$
(I^+ + II^+ + III^+) 
= \frac 1{2^3}\sum_{\sig_1, \sig_2=\pm 1}
\(
\nu  - \sigma_1(\mu_1+\rho_1)
-\sigma_2(\mu_2+\rho_2) + (\rho_1+\rho_2)
\) 
\phi_{\mu +\sig, l+1},   
$$
$$
(I^- + II^- + III^-) 
= \frac 1{2^3} \sum_{\sig_1, \sig_2=\pm 1}
\(\nu -
( \sigma_1(\mu_1+\rho_1) 
-\sigma_2(\mu_2+\rho_2) - (\rho_1+\rho_2) 
\) 
\phi_{\mu +\sig, l-1}.
$$
This finished the proof.
\end{proof}

\begin{rema+} It is remarkable that 
  in the formula the line bundle parameter 
  $l$ disappear
 due to the cancellation of $2l$ in the sum  $II+ III$. 
 When $\fg= \fsu(r+b, r)$, $r >1$, $\fk=\fs(u(r+b) +\fu(r))$, 
 the action of $\pi$ on $I(\nu)$ has been
 studied in details in  \cite{HT}. In this case
 the parameter $l$ indeed does not appear in affine term
 $\nu+\sigma_1(\mu_1+\rho_1)
+ \sigma_2(\mu_2+\rho_2) -\rho_1-\rho_2)$ for the action; 
 also the coefficients $-A^{\pm, \pm}$
 in \cite[Lemma 4.1]{HT}
 can be
formulated,  writing $\mu_1 =(\mu_1+\rho_1) -\rho_1,
\mu_2 =(\mu_2+\rho_2) -\rho_2$,  as
 $$
 \nu \pm  (\mu_1 +\rho_1)\pm (\mu_2+\rho_2) -\rho_1-\rho_2, 
 $$
 with our $\nu$ being their $-a=-(\alpha+\beta)$,
 $\mu_1=m_1+m_2$, $\mu_2=n_1+n_2$,
 $\rho_1=p-1$, $\rho_2=q-1$, and  our $l$  their $m_1-m_2=n_2-n_1$;
 see further \cite[(4.10)-(4.11)]{HT}.
 The Weyl group symmetry is again manifest here.
\end{rema+}

\section{Reductions points,   complementary and composition series}

\subsection{Reduction points and finite dimensional 
  subrepresentations}
We study now the existence of intertwining operators between
representations $I(\nu)$ and $I(\nu')$, and
we find certain finite dimensional 
representations at the reduction point of
$I(\nu)$.

\begin{theo+}
  \label{reductn}
   Let   $\fg   \ne \fsp(r, \mathbb R)$
    be a simple Hermitian Lie algebra of rank $r\ge 2$.
\begin{enumerate}  
  \item    There exists an intertwining operator between 
    the induced irreducible representations $I(\nu)$ and $I(\nu')$ 
  if and only if $\nu=\nu'$ or $\nu+\nu'=2\rho_{\fg}$.
\item    $I(\nu)$ is reducible if  and only if $\nu$ is an even
  integer, $\nu\ge 2\rho_1 +2$ or
  $\nu\le 2\rho_2 -2$.
    Moreover
  at the  point $\nu=-2k$, $k\ge 1$,  the
  the symmetric tensor product $S^k(\fg^{\mathbb C})$
  is  realized as (reducible, generally) finite-dimensional subrepresentation of $I(\nu)$
  via  
  \begin{equation}
    \label{symm-tensor-sub}
    T:  S^k(\fg^{\mathbb C})
    \to I(\nu), 
  X\mapsto f(g) =(\otimes^k \ad (g) (E_0), X),  g\in G,    
  \end{equation}
  where $E_0= \xi_{ie}-iH_0\in \fn_2$ is the  basis vector
  of the center $\nu_2$ of the nilpotent algebra $\fn=\fn_1+\fn_2$
  and $  (\cdot, \cdot)$ is the Killing form in 
  $ \fg^{\mathbb C}$ extended to   $  S^k(\fg^{\mathbb C})$.
  \end{enumerate}
\end{theo+}
\begin{proof}  The first part and the second part on reductions points
  are done similarly as in \cite{HT, Sahi-crelle,  Z-matann}.
  Now let  $\nu=-2k$ be an negative even integer.
  We prove that $T$ above is an intertwining operator from   $S^k(\fg^{\mathbb C})$
  into $I(\nu) , \nu=-2k$.
 The functions $f=f_X$ transform under $P=MAN$ as
$$
f(gm) =(\otimes^k \ad (gm) (E_0), X) 
=(\otimes^k \ad (g)\ad (m) (E_0), X) 
=(\otimes^k \ad (g), X) =f(g), \quad m\in M,
$$
since $M$ centralizes $E_0$,
$f(gn)=f(g)$, $n\in N$,  as $E_0$ is in the center of $N$,
and
$$
f(ge^{t\xi}) 
=(\otimes^k \ad (g)\ad (e^{t\xi}) (E_0), X) 
=(\otimes^k (e^{2t}\ad (g)(E_0) , X) =
e^{2kt}f(g)
=e^{-\nu t}f(g), \quad m\in N, 
$$
since    $\Ad (\xi)E_0=2E_0, \, \ad (e^{t\xi}) (E_0) =e^{2t} E_0$.

Thus $f\in I(\nu)$. The intertwining property of $T$
is obvious by its definition. This completes the proof.
\end{proof}

\subsection{Complementary series}
We  determine the complementary series i.e, that case when $\nu$
is real and the whole module
$(\fg^{\mathbb C}, K)$-module $I(\nu)$ is unitary and irreducible.

\begin{theo+} The complementary series $I(\nu)$ appears precisely in the range
  $\nu=\rho_{\fg}  +\delta$, $|\delta|<\delta_0$, 
  $$
 \delta_0=\begin{cases}
1+b, &    \fg=\fsu(r+b, r)\\
   3, & \fg=\fso^\ast(2r)\\
    n-3, & \fg=\fso(2, n), \quad n>4\\
    3, & \fg=\fe_{6(-14)}\\
            5, & \fg=\fe_{7(-25)}
    \end{cases}.
    $$
\end{theo+}  

\begin{proof} The abstract arguments in  determining
  the complementary series here is the same as in \cite{HT, Sahi-crelle, Z-matann}, so we
  will only present some brief computations.
  Let $\nu$ be real. Suppose
  the $(\fg^{\mathbb C}, K)$-module $I(\nu)$ is irreducible 
  with invariant Hermitian inner product  $\langle \cdot, \cdot \rangle_\nu$.
  By Schur lemma we have
  $\langle f, f \rangle_\nu= S(\nu, \mu, l)
  \Vert f\Vert^2$
  for all $f\in W_{\mu, l}$, where
$\Vert f\Vert^2$   is the norm square in $L^2(K/L)$
  and $S(\nu, \mu, l)$ is the Schur proportionality constant.
  Then
  $\pi_\nu(\xi)$ is skew symmetric and in particular
  $$
  \langle  \pi_\nu(\xi) 
  \phi_{\mu, l},
  \phi_{\mu+\sig, l+1}
\rangle_\nu
  =-\langle  
  \phi_{\mu, l}, \pi_\nu(\xi)\phi_{\mu+\sig, l+1 }\rangle_\nu
  $$
  Write the expansion of
$  \pi_\nu(\xi) 
  \phi_{\mu, l}$  in Theorem \ref{main} as 
$$
  \pi_\nu(\xi) 
  \phi_{\mu, l}
  =\sum_{\sig_1, \sig_2=\pm 1} A(\nu, \mu, l; \mu +\sigma, l+1) 
  \phi_{\mu +\sig, l+1} +
  \sum_{\sig_1, \sig_2=\pm 1} A(\nu, \mu, l; \mu +\sigma, l-1)  
  \phi_{\mu +\sig, l-1} 
$$
Thus the invariance of of the Hermitian form  above becomes
\begin{equation*}
   A(\nu, \mu, l; \mu +\sigma, l+1) 
  S(\mu +\sig, l+1) 
  =-    A(\nu, \mu+\sigma, l+1; \mu, l) 
 S(\mu +\sig, l), \,
 \end{equation*}
 \begin{equation*}
 A(\nu, \mu, l; \mu +\sigma, l-1) 
  S(\mu +\sig, l-1) 
  =-    A(\nu, \mu+\sigma, l-1; \mu, l)  
S(\mu, l).
    \end{equation*}
    That $I(\nu)$ is unitary and irreducible is
     equivalent to all Schur
     proportionality constants    $S(\mu, l)$
     are positive. It implies
     that $   A(\nu, \mu, l; \mu +\sigma, l+1) $
     and  $A(\nu, \mu+\sigma, l+1; \mu, l)$ have opposite signs, as well as
$        A(\nu, \mu, l; \mu +\sigma, l-1) $
and   $A(\nu, \mu+\sigma, l-1; \mu, l) $. This  determines
the range of $\nu$, given by the condition
$$
|\delta| < \text{min}\{\rho, \rho-2\rho_2, 2\rho_1 -\rho +2\}.
$$
By  case by case computations we get the range as claimed.
\end{proof}  

\begin{rema+}
  The complementary series for $SU(p, q)$
    has been found before in \cite[(ii)(b), p.~49; 5.2, p.~69]{KS},
    \cite[4.4]{HT}.
\end{rema+}

\subsection{Composition series and unitarizable subrepresentations}
The composition series for $I(\nu)$ at
reducible points $\nu$ is a bit involved. We shall only
determine the unitary subrepresentations
at some typical reduction points. The proof
of the following result is
done by examining the signs of
the coefficients $\nu +\sigma_1(\mu_1+\rho_1)
+\sigma_2(\mu_2+\rho_2) -\rho_1-\rho_2$
in Theorem   \ref{main}.

\begin{theo+}   Suppose $\nu\le 2\rho_2-2$ is an even integer.
  Then there are two unitarizable subrepresentations
  $S^{\pm}(\nu)\subset I(\nu)$ consisting of
  the  $K$-types
  $$
  S^{+}(\nu) =\sum_{(l, \mu), \mu_1-\nu_2 \ge -\nu+2\rho_1} W_{\mu,
    l},\quad
    S^{-}(\nu) =\sum_{(l, \mu), \mu_1-\nu_2 \le -\nu+2\rho_2} W_{\mu, l}.
  $$
\end{theo+}

\section{ The cases of
 $\fg=\mathfrak{su}(d, 1)$  and
  $\fg=\mathfrak{sp}(r, \mathbb R)$
 }

We treat
the remaining
cases when $\fg=\mathfrak{sp}(r, \mathbb R), \mathfrak{su}(d, 1) 
$ where the  space $K/L_0$ is the complex projective space $\mathbb P^{r-1}, \mathbb
P^{d-1}$.

\subsection{$\fg = \mathfrak{su}(d, 1)$
}
\label{sect-5-1}

This case is already treated in \cite{JW}
by using rather explicit computations 
of differentiating hypergeometric functions. We shall give  somewhat easier 
proof of their results by using our method above; this
avoids explicit computations involving special functions
and gives conceptual expression for the action of $\pi_\nu(\xi)$ 
in  terms of Harish-Chandra $c$-functions.

The Cartan decomposition
is $\fg =\fk +\fp=\fu(d) +\mathbb C^d$,
with  $\fk^{\mathbb C}=\mathfrak{gl}(d)$.  
The Jordan triple $V=\mathbb C^d$ with
 $\{ x, y, z \} = D(x, y)z= \langle x, y\rangle z
+\langle z, y\rangle x$. We fix the  tripotent $e=e_1$,
 the standard basis vector of 
 $\mathbb C^d$ and $\xi:=\xi_{e}$. The half-sum $\rho_\fg$
 is $\rho_\fg=d$.
 The space $S=K/L$ is the sphere $S$ in $\mathbb C^d$
 with $L$  the isotropic subgroup of $e\in S$,
 and $S_1=K/L_0$
 the projective space $\mathbb P(\mathbb C^d)$
 with $S\to S_1$ as a circle bundle over
 $\mathbb P(\mathbb C^d)$ by the defining map
 $z\mapsto [z]$.
 The tangent space of $S_1$ is realized as
 $T_{[e]}^{(1, 0)}S_1=
 \{  D(v,  e); v\in V_1= \mathbb C^{d-1}\}$,
 $T_{[e]}^{(0, 1)}S_1=
\{  D(e,  v); v\in V_1\}$.
 We fix an 
 $\fsu(2)$-subalgebra in $\fk$ and
 will perform the differentiation along
 tangent vectors $E^{\pm}$ in the complex totally geodesic
submanifold $\mathbb P^1\subset S_1$:
 $$
 E^+ = D(e_2, e), \, E^- = D(e, e_2), E=E^+ -E^-\in T_{[e]}(K/L_0).
 $$
Put $ H=[E^+, E^-]=D(e_2, e_2)-D(e, 
e)$.
 The element $iE$ generates a Cartan subalgebra
 for the  non-compact dual $(\fk^*_1, \fl_1)$, and
 positive roots are $\{2\alpha_1, \alpha_1\}$ with $\alpha_1$
 the dual element of $iE$, and the half-sum is $\rho_{\fk^*_1} = d-1$.

The decomposition of $L^2(K/L)$ is well-known,
  $$
  L^2(K/L) 
  =\sum_{m\ge |l|, m=l \,\text{mod}\, 2} W_{m, l}. 
  $$
 Each space 
  $W_{m, l}$
is generated by  ${\overline z_1}^{p}
 { z_2}^{q}$. In terms of
 the previous notation
 $$p=\frac{m+l}2, q=\frac{m-l}2, \quad m \ge |l|, m=l \quad \text{mod}\quad 2.
 $$

Now the coefficients in the expansion of  $\pi_\nu(\xi) 
\phi_{m, l}$
can be written as
$\(\nu + \pm (m +d-1) + c\) $ for some constants
$c$ as in Theorem \ref{main}. We write them explicitly.

 \begin{theo+}
  \label{main-2}
  \begin{enumerate}
  \item
    The action of $   \pi_\nu(\xi)$ on 
$\phi_{m, l}$ is given by 
\begin{equation*}
  \begin{split}
2^2 \pi_\nu(\xi) 
\phi_{m, l}
&\quad 
=\(\nu +m +l\) 
c_{m, l}(m+1, l+1) \phi_{m +1 , l+1}
+
\(\nu -m-2d +2 +l\) 
c_{m, l}(m+ 1, l+1) \phi_{m -1 , l+1}
\\
&\quad +
\(\nu + m  -l \) 
c_{m, l}(m+1, l-1) 
\phi_{m +1 , l-1}
+\(\nu - m  -l -2d+2 \) 
c_{m, l}(m-1, l-1) 
\phi_{m -1 , l-1}.     
  \end{split}
\end{equation*}
  \item The complementary series is in the range
    $\nu=\rho_{\fg} + \delta$, $|\delta| < d=\rho_{\fg}$.
\end{enumerate}
\end{theo+}

\begin{proof} We follow the computations of  $
  \pi_\nu(\xi) 
  \phi_{m, l}
  $ in the proof of Theorem \ref{main}
  and indicate the necessary changes. 
  We find first the coefficient of 
  $\phi_{m, l}$ and the other coefficients 
  will be found by general arguments. 
We  have 
 $ \pi_\nu(\xi) 
  \phi_{m, l}= I+ II + III$, $I=I^+ + I^-$, 
$$I^+=  \frac{\nu}{2} \phi_{1, 1}(k) 
\phi_{m, l}(k) 
$$
and the spherical polynomial $\phi_{1, 1}(k)$
is 
$$
\phi_{1, 1}(k)= 
\phi_{1, 1}(z) 
=\langle k e, e\rangle = z_1, \quad z=k e\in S. 
$$
as function on the sphere $S=K/L$. For $k=\exp(tE)$ we have 
$\phi_{1, 1}(k)= \cos x$ and its complexification is 
$\cosh x=\frac 12(e^x + e^{-x})$. The expansion 
(\ref{elem-phi}) now has coefficient $\frac 12 $. 
Thus the leading term in $I^+$ is 
$$
\frac{\nu} 4 c_{m, l}(m+1, l+1) \phi_{m+1, l+1}
$$
where $c_{m, l}(m+1, l+1)$ is  a quotient of two Harish-Chandra $c$-functions.

The  term $II^+$ in (\ref{II-upd-1}) in the present case becomes
$
II^+ = l \phi_{1, 1} \phi_{m, l}
$
and has the leading term 
 $$
\frac{l} 2 c_{m, l}(m+1, l+1) \phi_{m+1, l+1}.
$$

The third term $III^+$ is treated similarly as in Lemma \ref{diff}, and 
(\ref{III-key}) gives that the leading term of 
$III^+$
is
$\frac{m-l} 4 c_{m, l}(m+1, l+1) \phi_{m+1, l+1}.
$
Altogether terms involving
$\phi_{m+1, l+1}$ in $\pi_\nu(\xi)\phi_{m, l}$
are
$$
\frac{\nu +2l + m-l}{4}
c_{m, l}(m+1, l+1) \phi_{m+1, l+1}=
\frac{\nu +m+l}{4}
c_{m, l}(m+1, l+1) \phi_{m+1, l+1}
$$
with
$$
\frac{\nu +m+l}{4}
=\frac{\nu + (m+\rho_{\fk_1^*})
  -\rho_{\fk_1^\ast} +l}{4}
$$
in terms of the Weyl group invariant parameter $m+\rho_{\fk_1^*}$.
Observe again  that  
the highest weight with respect to $E$ is  $m\alpha_1$ 
so the Weyl group symmetric is with respect  
to $(m +\rho_{\fk_1^\ast})\alpha_1
$.
The coefficient of
$\phi_{m+1, l+1}$  involving $\nu$ is 
$$\nu + m +l=\nu + (m+ \rho_{\fk_1^\ast}) -\rho_{\fk_1^\ast} +l. 
$$
Thus $\phi_{m-1, l+1}$  has coefficient
$$\nu - (m+ \rho_{\fk_1^\ast}) -\rho_{\fk_1^\ast} +l
=\nu - m -2\rho_{\fk_1^\ast} +l =
\nu - m -2d+2 +l.
$$
The other two coefficient 
$\phi_{m\pm 1, l-1}$ is found 
by the unitarity of $\pi_\nu(\xi)$
at $\nu=\rho_\fg +i x, x\in \mathbb R$
and by the Weyl group symmetry.

The rest is obtained as in \cite{JW}.
\end{proof}

\begin{rema+}  The statement of \cite[Theorem 4.1]{JW}
  for the group $SU(n, 1)$  is (their $d=2$)
  \begin{equation*}
    \label{eq:1}
    \begin{split}
  \pi_\nu(H)   e_{m, l}
 & = a_{m+1, l+1}
  (\nu +m+l) 
  e_{m+1, l+1}
  + a_{m-1, l+1} (\nu -m+l-2n+2) e_{m-1, l+1} \\
  &
  + a_{m+1, l-1} (\nu  +m-l) e_{m+1, l-1}
  + a_{m-1, l+1} (\nu -m-l -2n+2) e_{m-1, l-1}
  \end{split}
  \end{equation*}
    where $a_{m\pm 1, l\pm1}$ are certain positive 
    constants independent of $\nu$. The important
  coefficients are   $\nu +m+l$
  for $  e_{m+1, l+1}$ and
  $\nu  +m-l$ for $e_{m+1, l-1}$,  which have
  rather simple form, and the rest
  is obtained by Weyl group symmetry.
    This 
    coincides then with our formula. 
    As mentioned above in Remark \ref{rema+-}
    it is enough to determine the leading 
    coefficient $a_{m+1, l+1} (\nu +m+l)$. 
\end{rema+}

\subsection{  $\fg = \mathfrak{sp}(r, \mathbb R)$, $r\ge 2$
}

The Jordan triple here is $V=
 M_r^s=\{v\in M_r(\mathbb C); v=v^t\}$ of 
complex matrices with the triple product $D(u, v)w=u v^\ast w + w v^\ast u$. This
case is rather special
and we provide all details.
The normalization of Euclidean norm in $\fp^+$ is as before
with  minimal tripotents having  norm $1$. 
The group $K=U(r)$
acts on $V$
and by $A\in K: Z\to AZA^t$. To avoid confusion with 
various realizations we recall that
all Lie algebra elements are realized via Jordan triple
as in Section 1, in particular  Lie algebra elements
of $\fk^{\mathbb C}=\fgl(r, \mathbb C)$
appear as $D(u, v): w\to D(u, v)w$; they are identified with usual matrices $uv^*$
if we still want matrix realizations.

We fix the minimal tripotent
the diagonal matrix $e=\diag(1, 0, \cdots, 0)\in V$ and $\xi=\xi_e\in \fp$
 in    (\ref{xi-H0}).
The functional $\rho_\fg$ is now $\rho_{\fg} =r$.
A subtle point here is that the group $L=M\cap K\subset K$ is 
$$
L 
=\mathbb Z_2\times U(r-1)=\{h=\diag(h_0, h_1); h_0=\pm 1, h\in U(r-1)\}.
$$
Let 
$$
L_0=U(1)\times U(r-1)=\{h=\diag(h_0, h_1); h_0\in U(1), h\in U(r-1)\}.
$$
Thus 
$S=K/L\subset V$ is the real projective space 
$\mathbb P(\mathbb R^{2r})\subset V$ and $S_1=K/L_0$
is again the complex projective space 
$ \mathbb P(\mathbb C^r)\subset \mathbb P(V)$ of rank one realized
in the projectivization of $V$.

The Peirce 
decomposition $V=V_2 +V_1 + V_0$ with respect to $e$ is the block $(1 + r-1)\times
(1\times r-1)$-partition of $V$. We fix 
$$v=\begin{bmatrix}0 &1 \\1&0\end{bmatrix}\in
V_1,\, 
w=\begin{bmatrix}0 &0 \\0&1\end{bmatrix}\in   V_0,\, 
$$ 
$v$ is a tripotent of {\it rank  two}
and $w$ of {\it rank one}. 
The Cartan decomposition of $\fk$
is $\fk = \fu(r-1) + \fq$ where 
elements in $\fq=\{D(u, e)-D(e, u); u\in V_1\}=\mathbb  C^{r-1}$  as real spaces.

We fix also the following  $\fsl(2)$ elements,
\begin{equation}
  \label{sp-H-12}
E = D(v, e) - D(e, v) =E^+ -E^-
\in \fq, \quad H= [E^+, E^-]=D(v, v)-2D(e, e) 
\end{equation}
with 
$[H, E^+]= 2 E^+$.
As matrices $D(e, v)= E_{12} \in \fgl(r, \mathbb C).$
(Note the difference between this case and
\ref{H-12}), and this has been studied
in greater details in \cite{N}.) The 
symmetric pairs
$(\fk^*, \fl_0)=\fu(r-1, 1), \fu(r-1) +\fu(1))$,
$(\fk^*_1, \fl_1)=(\fsu(r-1, 1), \fs(\fu(r-1) +\fu(1))  )$,
and the roots of
$(\fk^*_1, iE)$ are $2$ and $1$ with 
$\rho_{\fk_1*\ast}=r-1$, and $\rho_\fg= 1 +\rho_{\fk_1^\ast}$.

Lemma \ref{C-H} in the present case is 
  $$
  L^2(K/L) 
  =\sum_{m\ge |l|}  W_{2m, 2l}, 
  $$
where each space 
  $W_{2m, 2l}$ is generated by
 ${\overline  z_1}^{(m+l)}
  { z_2}^{(m-l)}$ and  contains  the spherical polynomial
  $\phi_{2m, 2l} $.

\begin{theo+} 
  \label{xi-action}
\begin{enumerate}
\item 
The action of $   \pi_\nu(\xi)$ on 
$\phi_{2m, 2l}$ is given by
\begin{equation*}
  \begin{split}
2^3  \pi_\nu(\xi) 
\phi_{2m, 2l}
&=\(\nu   +2m\) 
c^{m+1, l+1}_{m, l}
\phi_{2m +2\sig , 2l+2} + \(\nu   -2m-2r+2\) 
c^{m-1, l+1}_{m, l}
\phi_{2m -2 , 2l+2} 
\\
&+ 
\(\nu +  2m\) c^{m+1, l-1}_{m, l}
\phi_{2m +2 , 2l-21}
+\(\nu -  2m -2r +2\) 
c^{m-1, l-1}_{m, l}
\phi_{2m -2 , 2l-21}
  \end{split}
\end{equation*}
\item
  There is no complementary series
  in the family $I(\nu)$ if $r=2$.
  If $r>2$ the complementary series are in the range
  $\nu=r+\delta$, $|\delta| < r-2$.
\end{enumerate}
\end{theo+}

\begin{proof} Following the earlier computation
  in Section 4 we  have 
$$I^+=  \frac{\nu}{2} \langle k e, e\rangle
\phi_{2m, 2l}(k). 
$$
The matrix coefficient
$\langle k e, e\rangle
=\phi_{2, 2}(k)$, 
and its restriction on the torus
$\exp(\mathbb RE) e$ is
\begin{equation}
  \label{exp-xE-e}
\exp(xE) e=
\begin{bmatrix} \cos^2 x &\cos x \sin x \\
\cos x \sin x& \sin^2 x\end{bmatrix}
=\cos^2 x e +\cos x \sin x v + \sin^2x w.   
\end{equation}
Thus
$\phi_{2, 2}(k)= (\cos x)^2$ and it has
an expansion
$$
\phi_{2, 2}(k)= (\cos x)^2=\frac{1}4 (e^{2ix} + 2 + e^{-2ix}).
$$
Thus the leading term in $I^+$ is 
$$
\frac{\nu} 8 c_{m, l}(m+1, l+1) \phi_{2m+2, 2l+2}
$$
where $c_{m, l}(m+1, l+1)$ is the quotient Harish-Chandra $c$-functions 
for  $\phi_{2m+2, 2l+2}$ and $\phi_{2m, 2l}$.

In the second term $II$ we l have $H_0e =2ie$, 
and then $II^+ = l \phi_{2, 2} \phi_{2m, 2l}$, which
has the leading term 
 $$
\frac{l} 4 c_{m, l}(m+1, l+1) \phi_{2m+2, 2l+2}.
$$

The vector field $X(k)=D(e, P_1(k^{-1}e))$ in      (\ref{vect-X})
has restriction
$$
X(\exp(xE)
)=D(e, P_1(\exp(-xE) e)) = -\cos x\sin x D(e, v), 
$$
 by (\ref{exp-xE-e}). 
 The leading term of the expansion  $-\cos x\sin x E^-\phi_{2m, 2l}$ is obtained,
 using the proof of Lemma \ref{A1} 
 in Appendix \ref{app-a}, as
 $$
 -\cos x\sin x E^-\phi_{2m, 2l} (\exp(xE))
 =\frac {2m-2l}8 e^{i2mx} + L.O.T..
 $$
 Altogether we find
 \begin{equation*}
   \label{eq:3}
   \begin{split}
 I^+ + II^+ + III^+
&= \frac{\nu+ 2l +2m-2l} 8 c_{m, l}(m+1, l+1) \phi_{2m+2, 2l+2}  + L.O.T 
\\
&= \frac{\nu+ 2m} 8 c_{m, l}(m+1, l+1) \phi_{2m+2, 2l+2}  + L.O.T..     
   \end{split}
 \end{equation*}
The rest is done as  in the proof of Theorem \ref{main} above.
\end{proof}

\begin{rema+} We note  the trivial representation corresponds to
  $\nu=r+\delta$, $\delta =r$,   
  whereas the complementary serie range is $|\delta| < r-2$, so the 
  gap  is $2$ between the end of complementary series
  and the trivial representation, similar to  the rank one case $Sp(n, 1)$
  (and there is no gap for $SU(n, 1)$).
\end{rema+}  

\begin{rema+}
Even integers  $\nu=-2k$ are reduction points for $I(\nu)$ 
in both cases above,  $\fg=\fsu(d, 1), \fsp(r, \mathbb R)$. 
The  map  (\ref{symm-tensor-sub}) 
realizes $S^k(\fg^{\mathbb C})$ as 
a subrepresentations of $I(\nu)$. 
\end{rema+}

There are some interesting non-spherical unitary
principal series
of $Sp(n, \mathbb R)$ and it might be interesting
to study them in details.

\appendix 
\section{Recursion formulas for differentiations of spherical 
  polynomials}
\label{app-a}
Let
$$SU(2)=\left\{g=
\begin{bmatrix} a&b \\
  -\bar b & \bar a
  \end{bmatrix}; |a|^2+|b|^2=1\right\}
  $$
  and fix the following  Lie algebra $sl(2, \mathbb C)$ elements, 
$$
H=\begin{bmatrix} 1&0 \\
  0 & -1 
\end{bmatrix},\, 
E^+ =
\begin{bmatrix} 0&1 \\
  0 & 0 
\end{bmatrix}, 
E^- =
\begin{bmatrix} 0&0 \\
  1 & 0 
  \end{bmatrix},   E=E^+ + E^-\in i\fsu(2). 
  $$
Then the   $\fsl(2)$-algebras
    $\fsl(2)=\mathbb C H_j + \mathbb C E_j^+ +
    \mathbb C E_j^-, j=1, 2,$
    defined in 
    (\ref{H-12}) are isomorphic
to the present $\fsl(2)$ via the identification
$$
H_j\longleftrightarrow H,   E_j^{\pm}\longleftrightarrow 
E^{\pm};
$$
as well as the identification of the compact $\fsu(2)$-real forms
$$
\fsu(2)=\mathbb R iH_j  +  \mathbb R E_j +\mathbb R(i(E_j^+ +E^-_j))
\longleftrightarrow  \mathbb R iH  +  \mathbb R E +\mathbb R(i(E^+ +E^-)).
$$
Let 
 $$
U(1)=\left\{ u_\theta=
\begin{bmatrix} e^{i\theta}
    &0 \\
 0 & e^{-i\theta}
  \end{bmatrix} =e^{i\theta H}
  \right\}. 
  $$
  Recall  $\chi_l(H_j)=-l$ and the transformation rule
  (\ref{chi-H-j})  of   $\phi_{\mu, l}$ under $\exp(itH_j)$.
  Accordingly we let $\chi_l(H)=-l$, and the spherical polynomials $\phi_{m, l}$ in the present
  case satisfy
$\phi_{m, l}(e^{i\theta H}g) =\phi_{m, l}(ge^{i\theta H}) =e^{-il\theta}\phi_{m, l}(g).
$ 
 Consider   the representation of $SL(2, \mathbb C)$ on the symmetric tensor 
  $S^m:=\bigodot^m \mathbb C^2$
  of the defining representation $\mathbb C^2$,
and write  the action simply as $g\to gv$, $v\in  S^m \mathbb
C^2$,
as well as the Lie algebra action.
  Let $-m\le l \le m$, $m=l$ mod $2$.
    The $l$-spherical polynomial
    is given by the matrix coefficient
    \begin{equation}
      \label{phi-sl2}
  \phi_{m, l} (g) 
  =\binom{m}k 
  \langle  g (e_1^ke_2^{m-k}), e_1^ke_2^{m-k}
  \rangle, 
  0\le k\le m, \,  {l=m-2k},       
    \end{equation}
    where  the tensor $e_1^ke_2^{m-k} = e_1^k\otimes e_2^{m-k}$
    is as usual
  viewed as polynomial on the dual space of $\mathbb C^2$.
This is  verified by
\begin{equation}
    \begin{aligned}  
      \phi_{m, l} (gu_\theta) =
      \phi_{m, l} (g e^{i\theta H})
      &  =
  \langle   g e^{i\theta H}(e_1^ke_2^{m-k}), e_1^ke_2^{m-k}
\rangle
=e^{i(2k-m)\theta}\phi_{m, l} (g)
\\
&=e^{-il\theta}\phi_{m, l} (g)
=\chi_l(u_\theta)
\phi_{m, l} (g),
\end{aligned}
\end{equation}
also $  \phi_{m, l} (u_\theta g) =
\chi_l(u_\theta)\phi_{m, l} (g)
$ along with the normalization $\phi_{m, l}(I) = \binom{m}{k}
\Vert e_1^ke_2^{m-k}\Vert^2 =1$.
  For our purpose in Section \ref{Liealg-act} it is more convenient
  to consider the spherical polynomial
\begin{equation}
  \label{psi-sl2}
  \psi_{m, l}(g)
  =\frac{2^m}{\binom{m}k}
  \phi_{m, l} (g)= 2^m   \langle  g (e_1^ke_2^{m-k}), e_1^ke_2^{m-k}
  \rangle, \,  {l=m-2k},       
  \end{equation}
    which has the normalization the leading
    term of $\psi_{m, l}(\exp tE)$ being $e^{imt}$, i.e.
    $$
    \psi_{m, l}(\exp tE) = e^{imt} + L.O.T.
    $$
where L.O.T. is a trigonometric polynomial of lower order.
The following lemma is used in the proof of Lemma \ref{diff}. Actually we need
only to find the leading term of the trigonometric polynomials
$\langle  g e_1, e_2\rangle  (E^-\phi_{m, l}) (g)$ and
$  \langle  g e_2, e_1\rangle  (E^+\phi_{m, l}) (g) $
for $g=\exp(tE)$, which is  elementary.
\begin{lemm+} \label{A1} 
The following recurrence formulas hold,
\begin{equation}
  \label{E+phi}
 \langle  g e_1, e_2\rangle  (E^-\phi_{m, l}) (g) 
  =
 \frac 14 
  \frac{(m - l)(m+l+2)}
  {m+1}
  \left(
    \phi_{m+1, l+1}
- \phi_{m-1, l+1}\right),
\end{equation}
  \begin{equation}
  \label{E-phi}
  \langle  g e_2, e_1\rangle  (E^+\phi_{m, l}) (g) 
  =
  \frac{1}4 
  \frac{(m+l)(m-l+2) 
  } {m+1}
\left( \phi_{m+1, l-1} -  \phi_{m-1, l-1}\right).
  \end{equation}
When restricted to  
$g=\exp(xE)=\begin{bmatrix} \cos x& \sin x\\
  -\sin x& \cos x\\
  \end{bmatrix}
  $ and written in terms of $\psi_{m, l}$  they are
\begin{equation} 
  \label{E+psi-expl}
- \sin x  (E^-\psi_{m, l}) (g) 
=\frac 14 (m-l) 
\psi_{m+1, l+1}  -\frac 14 \frac{(m+l+2)(m-l)^2}{m(m+1)} 
\psi_{m-1, l+1},
\end{equation}
\begin{equation}
  \label{E-psi-expl}
  \sin x (E^+\psi_{m, l})(g) =
  \frac 14 (m+l)  \psi_{m+1, l-1}(g) 
  - \frac 14 \frac{(m-l+2)(m+l)}{m(m+1}
   \psi_{m-1, l-1}(g).
 \end{equation}
\end{lemm+}

\begin{proof} First we find the weight
  of $f(g)=\langle  g e_2, e_1\rangle  (E^+\psi_{m, l}) (g)$
  under the regular left
  and right actions  of $\exp(i\theta H)$.
  We have   $E^+$ is of weight $2$ under $\ad(H)$.
    Also
  the matrix coefficient    $\langle g e_2, e_1\rangle $ transforms
  as
  $\langle  ge^{i\theta H} e_2, e_1\rangle =
e^{-i\theta } \langle  g 
e_2, e_1\rangle  $,
thus it is of weight $-1$ under the right regular action
of $H$.
The character $\chi_l$ is defined by $\chi_l(H)=-l$
thus   $f(g)=\langle  g e_2, e_1\rangle  (E^+\psi_{m, l}) (g)$
is of weight $-l +2-1 =-(l-1)=\chi_{l-1}(H)$
in the sense $f(g\exp(i\theta H)) =
\chi_{l-1}(\exp(i\theta H)) f(g)$.
Similarly $\langle  \exp({i\theta H}) ge_2, e_1\rangle =
e^{i\theta } \langle  g 
e_2, e_1\rangle  $
and the right differentiation by $E^+$
commutes with the left action. Thus
$f( \exp({i\theta H}) g)  = e^{-i(l-1)\theta}f(g)$, 
and $f(g)$   is  a linear combination of
$ \phi_{m+1, l+1}$ and $ \phi_{m-1, l+1}$ as it is the matrix
coefficient of  $\mathbb C^2\otimes S^m
  = S^{m+1}  \oplus
  S^{m-1} \mathbb C^2$.
  \begin{equation}
    \label{eq:A-B}
f(g)
=A \phi_{m+1, l-1} + B \phi_{m-1, l-1}.
  \end{equation}
    We have   $\phi_{m, l} (g)=\binom{m}k
  \langle  g (e_1^ke_2^{m-k}), e_1^ke_2^{m-k}
  \rangle
  $, and 
  $$(E^+\phi_{m, l}) (g)
  =\binom{m}k
  \langle  g (E^+ (e_1^ke_2^{m-k})), e_1^ke_2^{m-k}
  \rangle
  =\binom{m}k (m-k)
  \langle  g (e_1^{k+1}e_2^{m-k-1}),  e_1^{k}e_2^{m-k}
  \rangle,
  $$
  since $E^+e_1^ke_2^{m-k} =(m-k) e_1^{k+1}e_2^{m-k-1}$.
  Its product with
  $\langle  g e_2, e_1\rangle
$ is
  $$
  \langle  g e_2, e_1\rangle
  \langle  g (e_1^{k+1}e_2^{m-k-1}),  e_1^{k}e_2^{m-k}
  \rangle
  =  \langle  g (
  e_2\otimes  
  e_1^{k+1}e_2^{m-k-1}),   e_1\otimes e_1^ke_2^{m-k},
    \rangle.
  $$
Both $e_1\otimes e_1^ke_2^{m-k}$ and 
$e_2\otimes
e_1^{k+1}e_2^{m-k-1}$ are of weight $i(1+2k-m)$
under $H$,   the corresponding
weight vector in the space
$S^{m+1} $ respectively
$S^{m-1} $  is
$e_1^{k+1}e_2^{m-k}$ resp. 
$e_1^k e_2^{m-k-1}$ with matrix coefficient
$\langle  g (e_1^{k+1}e_2^{m-k}), 
e_1^{k+1}e_2^{m-k}\rangle$, 
$  \langle  g (e_1^k e_2^{m-k-1}),
e_1^{k}e_2^{m-k-1}\rangle $.  
In view of       (\ref{phi-sl2})
the formula  (\ref{E-phi}) becomes
\begin{equation}
  \begin{aligned}
&\quad      \binom{m}k (m-k) 
\langle  g (
e_2\otimes 
e_1^{k+1}e_2^{m-k-1}), e_1\otimes e_1^ke_2^{m-k},
\rangle 
 \\
 &=A \binom{m+1}{k+1} \langle  g (e_1^{k+1}e_2^{m-k}), 
  e_1^{k+1}e_2^{m-k}\rangle 
  +B \binom{m-1}{k}
  \langle  g (e_1^k e_2^{m-k-1}), 
   e_1^{k}e_2^{m-k-1}\rangle.
\end{aligned}  
\end{equation}
  Evaluating at  $g=I$  we get $B=-A$.
  Next we specify the equality  to
  the self adjoint element
  $$g=\begin{bmatrix} \ch x & \sh x\\
    \sh x & \ch x\\
  \end{bmatrix}
  =\begin{bmatrix} \ch \frac x2 & \sh \frac x2\\
    \sh \frac x2 & \ch \frac x2\\
  \end{bmatrix}^2 := h^2
    $$ and look for the coefficients
    of $e^{(m+1)x}$. We have
    $$
    \langle 
    g (
e_1\otimes e_1^ke_2^{m-k}), e_2\otimes 
e_1^{k+1}e_2^{m-k-1}\rangle
=   \langle 
    h (
e_1\otimes e_1^ke_2^{m-k}), h(e_2\otimes 
e_1^{k+1}e_2^{m-k-1})\rangle,
$$
and  its leading term 
 is
$$
\frac{e^{(m+1)x}} {2^{2(m+1)}}
\langle (e_1 + e_2)^{m+1},  (e_1+e_2)^{m+1}
\rangle
=\frac{e^{(m+1)x}}{2^{m+1}},
$$
and the LHS has leading term
$
\binom{m}k (m-k)
\frac{e^{(m+1)x}}{2^{m+1}}
.$
The term $e^{(m+1)x}$ appears only in the first summand
in the RHS which  has leading term
$A \binom{m+1}{k+1}\frac{e^{(m+1)x}}{2^{m+1}}.$
Thus
    $$
    A=\frac{
      \binom{m}k (m-k)}{\binom{m+1}{k+1}}
    =    \frac{(m-k)(k+1)}{m+1}
    = \frac{(m+l)(m-l+2)}{4(m+1)}.
    $$
    This proves (\ref{E-phi}), and then 
(\ref{E-psi}) by  using   (\ref{psi-sl2}) 

The  formulas (\ref{E+phi}) and
(\ref{E+psi})
are proved by the same 
 methods.

\end{proof}

\begin{rema+}
  In terms of $\psi_{m, l}$ they become 
\begin{equation}
  \label{E+psi}
\langle  g e_1, e_2\rangle  (E^-\psi_{m, l}) (g) 
=\frac 14 (m-l) 
\psi_{m+1, l+1}  -\frac 14 \frac{(m+l+2)(m+l)^2}{m(m+1)} 
\psi_{m-1, l+1}. 
\end{equation}
\begin{equation}
  \label{E-psi}
\langle  g e_2, e_1\rangle  (E^+\psi_{m, l}) (g) 
=\frac 14 (m+l) 
\psi_{m+1, l-1} - \frac 14 \frac{(m+2-l)(m+l)^2}{m(m+1)}
\psi_{m+1, l-1}; 
\end{equation}

The spherical polynomial
$\phi_{m, l}$
is special case of the
 spherical function 
 $\Phi_{\lambda, l}$ with $\lambda =-i(m +\rho)=-i(m +1)$
 in our case,
 and 
  $\Phi_{\lambda, l}$ is invariant with respect to the Weyl group
  action
  $\lambda\to -\lambda$. 
Namely, the pair of the coefficients
  $$
\pm  \frac{1}4
  \frac{(m-l))(m+l+2)  }
  {m+1}=
  \pm \frac{1}4
  \frac{((m+1)-(l+1))
    ((m+1)+(l+1))
  }
  {m+1}
$$ 
in the lemma  is invariant by the change $m+1\to -(m+1)$
and this symmetry is indeed  obvious here.  These formulas
are all classical trigonometric identities and
can be obtained by other methods.

\end{rema+}

\section{Table of Hermitian symmetric spaces $G/K$
  and their varieties of minimal rational tangents $K/L_0$.
  Duality relation for $(\dim(X), \genus(X))$ for $X=G/K, K/L_0$
}
\label{app-b}
\subsection{Tables}
We give a list of $G/K$
and the corresponding
projective
spaces  $S_1=\mathbb P(S) 
=
K/L_0=K_1/L_1$  as  compact Hermitian symmetric space; see \cite{He2, Loos-bsd, Hwang-Mok}.
The compact dual of
a noncompact Hermitian symmetric space $D$
is denoted by $D^\ast$. 

\begin{table}[!h]
      \begin{center}
\begin{tabular}
{
|l|l|l|l|
}
\hline
  $D=G/K$ & $G$ & $K$ & 
 $(a, b)$ 
\\
\hline
$I_{r+b, r}$
& 
              $ SU(r+b, r))$ &
                               $
                S(U(r+b)  \times U(r) )$ 
  &$ (2, b)$
   \\ 
  \hline
  $II_{2r}$
& 
                                    $ SO^\ast (4r)$
            &

                                               $ U(2r)$ 
&$    (4, 0) $
  \\ 
  \hline
    $II_{2r+1} $
& 
                                      $ SO^\ast (4r+2)$

            &$ U(2r+1)$ 
&
$ (4, 2)$
  \\
    \hline
  $III_{r} $
& 
              $ Sp(r, \mathbb R)$
            &  $ U(r)$ 
  &
$ (1, 0)$                                    
  \\ 
  \hline
  $IV_{n}, n>4. (r=2) $
& 
                   $SO(n, 2)$ &
                                $ SO(n)\times SO(2)$ 
&$ (n-2, 0)$                                  
  \\ 
\hline  
  $V  (r=2) $
& 
         $E_{6(-14)}$ &
                 $ Spin(10)\times SO(2)$ 
&$ (6, 4)$                                  
  \\ 
\hline  
  $VI (r=3) $
& 
       $E_{7(-25)}$ &
                  $E_6\times SO(2)$ 
                            &$ (8, 0)$
                                  
  \\ 
\hline  
\end{tabular}
\vskip0.20cm
 \caption
 {Non-compact Hermitian symmetric space $D=G/K$}
\label{tab:1}
      \end{center}
\end{table}

\begin{table}[!h]
      \begin{center}
\begin{tabular}
{
|l|l|l|
}
\hline
  $D=G/K$ &  $\mathbb P(S)=K/L_0=K_1/L_1$
    & $(a_1, b_1)$
\\
\hline
$I_{r+b, r}$
&          $I_{r+b-1}^*\times
                         I_{r-1}^*                                         
              $
     & {$ (0, r+b-2), (0, r-2)$}
    \\ 
  \hline 
  $II_{2r}$
&   $I^\ast_{2, 2r-2}  $
          &$ (2, 2r-4)$
  \\
  \hline 
  $II_{2r+1}$
&   $I^\ast_{2, 2r-1}  $
          &$ (2, 2r-3)$
       \\   
  \hline
  $III_{r} $
& 
            $I_{r-1}^*   $
    &$ (0, r-2)$
                                  
  \\ 
  \hline
  $IV_{n}, n>4 $
&    $IV_{n-2}^\ast $
    &$ (n-4, 0)$
                                  
  \\ 
\hline  
  $V $
& 
               $II_{5}^\ast $
                            &$ (4, 2)$
  \\ 
\hline  
  $VI $
& 
                         $V^\ast$ 
                                    &$ (6, 4)$
                                  
  \\ 
\hline  
\end{tabular}
\vskip0.20cm
 \caption
 {The   compact Hermitian symmetric space
   $\mathbb P(S)=K/L_0=K'/L'$. For type I domain $I_{r, r+b}$, $r\ge 2$,
   $\mathbb P(S)$ is a product $\mathbb P^{r-1}\times
\mathbb P^{r+b-1} $
   of projective spaces with the corresponding $(a_1, b_1)$
   is $(0, r+b-2), (0, r-2)$  for each factor.
}
\label{tab:2}
      \end{center}
\end{table}

\subsection{Duality between $(d, p)$
  and $(d_1, p_1)$ for $G/K$ and $K/L$}
Let  $d=\dim_{\mathbb C} D=r + \frac 12 a r(r-1) +r b$, $p 
  =2+a(r-1)+b$, the dimension and the genus of $D$. In terms
of Lie algebra actions they are
$$
(d, p) = (\tr \Ad (-iZ)|_{\fp^+}, 
\tr \Ad (D(e, e))|_{\fp^+})
$$
where $D(e, e)$ is the Harish-Chandra co-root of $\gamma_1$,
$\gamma_1(D(e, e))=2$.
Similarly let 
$$(d_1, p_1) =
(\dim(K_1^\ast/L), \text{genus}(K_1^\ast/L))$$
if $D\ne 
SU(r, r+b)/S(U(r)\times U( r+b))$. 
Put  $$
d' =\dim (\mathbb P^{r-1})= r-1, p'=r;  
d'' =\dim (\mathbb P^{r+b-1})= r+b-1, p''=r+b$$
when 
$D=SU(r, r+b)/S(U(r)\times U( r+b))$.

The
following duality  between the pairs 
$(\dim(D), \text{genus}(D))$
and $(\dim(K/L_0), \text{genus}(K/L_0))$
is mentioned in Lemma \ref{C-H} and might be of independent interest.
It can be proved by trace computations or by  case-by-case computations
of the tables above.
 \begin{lemm+}
        \begin{enumerate}
\item Let  $D$ be of rank $r\ge 2$
and is one of the domains $II, IV, V, VI$.  Then
$$     
      \frac{p}{d} + \frac{d_1}{p_1}  =2. 
          $$ 
\item Let  $D$ be of Type $I$  with $r\ge 2$.
Then            $$
          \frac{p}{d} + \frac{d'}{p'} +
           \frac{d''}{p''}     =2. 
          $$
\item  Let $D$ be  the Siegel domain $II$. Then
$$
\frac{p}{d} + 2\frac{d_1}{p_1} =2. 
$$
\end{enumerate}
 \end{lemm+}

\newpage

\def\cprime{$'$} \newcommand{\noopsort}[1]{} \newcommand{\printfirst}[2]{#1}
  \newcommand{\singleletter}[1]{#1} \newcommand{\switchargs}[2]{#2#1}
  \def\cprime{$'$} \def\cprime{$'$} \def\cprime{$'$} \def\cprime{$'$}
\providecommand{\bysame}{\leavevmode\hbox to3em{\hrulefill}\thinspace}
\providecommand{\MR}{\relax\ifhmode\unskip\space\fi MR }
\providecommand{\MRhref}[2]{%
  \href{http://www.ams.org/mathscinet-getitem?mr=#1}{#2}
}
\providecommand{\href}[2]{#2}


\begin{thebibliography}{1}






\bibitem{Barbasch}
  D. {Barbasch},
  \emph{
    The unitary spherical spectrum for split classical groups,
    } J. Inst. Math. Jussieu {\bf 9} (2010), no. 2, 265-356.

\bibitem{Frahm}J. Frahm, 
  \emph{    Conformally invariant differential operators on Heisenberg groups and minimal representations, }
  preprint, arXiv: 2012.05952v1.
    
\bibitem {He2}S.~Helgason, \emph{Groups and geometric analysis}, Academic
Press, New York, London, 1984.

\bibitem {He3}S.~Helgason, \emph{Geometric analysis on symmetric spaces},
Mathematical Surveys and Monographs, vol.~39, American Mathematical Society,
Providence, RI, 1994.
  

  
\bibitem 
{Howe}
R.~Howe, 
\emph{
  Some recent applications of induced representations,}
Group representations, ergodic theory, and mathematical physics: a 
tribute to George W. Mackey, Contemp. Math. {\bf 449} (2008), 
173-191.

\bibitem{HT}
  R.~Howe and E.~Tan,
  \emph{Homogeneous functions on light cones: The infinitesimal structure of some degenerate principal series representations},
  Bull. Amer. Math. Soc., {\bf 28} (1993), No. 1, 1-74.
    
\bibitem{Hwang-Mok}
  J.-M.~Hwang and N.~Mok, 
  \emph{
    Rigidity of irreducible Hermitian symmetric spaces of the 
    compact type under 
    K\"a{}hler deformation,}
  Invent. Math. {\bf 131} (1998), 393-418. 

  \bibitem{JW}
    K.~D.~Johnson and N.~R.~Wallach,
    \emph{
      Composition series and
  intertwining operators for the spherical principal series. {I}}, Trans. Amer.
  Math. Soc. \textbf{229} (1977), 137--173. 

\bibitem{Kn}
  A.~Knapp,
  \emph{
    Representation theory of semisimple groups.
An overview based on examples,  }
 Princeton Landmarks in Mathematics. Princeton University Press, Princeton, NJ, 2001.

  
\bibitem{KS} A.~Knapp and B.~Speh,
  \emph{
    Irreducible Unitary Representations of $SU(2, 2)$,
    }    J. Funct. Anal. \textbf{45}(1982), 41-73 .

\bibitem{KS-2} \bysame,
  \emph{
    The role of basic cases in classification: theorems about unitary representations applicable to SU(N,2),}
  Noncommutative harmonic analysis and Lie groups (Marseille, 1982), 119-160,
  Lecture Notes in Math., 1020, Springer,  1983.
    
\bibitem{Loos-bsd}O. Loos,
  \emph{Bounded Symmetric Domains and {Jordan} Pairs,}
University of California, Irvine, 1977.
    

\bibitem{Mok-geo-st}
  N.~Mok,
  \emph{
    Geometric structures on uniruled projective manifolds defined by their varieties of minimal rational tangents},
  Asterisque,
  {\bf 322} (2008), Volume II, 151-205.

\bibitem{N}{E. Neher}, 
  {\it Jordan Triple Systems by the Grid Approach,} Springer Lect. Notes in Math. {\bf 1280} (1987). 





  
\bibitem{OZ-duke}
  B.~\O{}rsted and G.~Zhang,
  \emph{  Generalized principal series representations and tube
    domains,}
  Duke Math. J. {\bf 78}
  (1995),  335-357.



\bibitem 
  {Sahi-crelle}
  S.~Sahi,
  \emph{ Jordan algebras and degenerate principal series},
  J. Reine Angew. Math. {\bf 462} (1995), 1-18.


\bibitem{Sch}H.~Schlichtkrull, 
  \emph{One-dimensional 
    K-types in finite dimensional 
    representations of semisimple 
    Lie groups: A generalization 
    of Helgason's theorem}, 
  Math. Scand. {\bf 54} (1984), 279-294.


\bibitem{Shi} N.~Shimeno, 
  \emph{
The Plancherel formula for spherical functions with a one-dimensional
K-type on a simply connected simple Lie group of Hermitian type,
  } 
J. Funct. Anal. {\bf 121} (1994), 330-388.

\bibitem{Upmeier}
  H. Upmeier,
  \emph{ Harald Jordan algebras and harmonic analysis on symmetric spaces,}
  Amer. J. Math. {\bf 108} (1986), no. 1, 1-25.

\bibitem{Vretare}
L.~Vretare,
\emph{
  Elementary spherical functions on symmetric spaces,}
  Math. Scand. {\bf 39} (1976),  343-358.


\bibitem{Z-matann}
  G.~Zhang,
  \emph{Jordan algebras and generalized principal series
    representations,}
  Math. Ann. {\bf 302} (1995),  773-786.

  \bibitem{Z-tams}
  \bysame,
  \emph{
    Some recurrence formulas for spherical polynomials on tube domains,
  }
     Trans. Amer. Math. Soc. 
 {\bf 347} (1995),  1725-1734.



\end{thebibliography}
\end{document}